\documentclass[a4paper,11pt]{article}
\usepackage[T1]{fontenc}
\usepackage[utf8]{inputenc}
\usepackage{lmodern}
\usepackage{amsmath,amssymb,amsthm,mathrsfs,mathtools}
\usepackage[margin=28mm]{geometry}
\usepackage{microtype}
\usepackage{needspace}
\usepackage{xcolor}
\usepackage[sort&compress,numbers]{natbib}
\usepackage{xurl}
\usepackage[hidelinks]{hyperref}
\newcommand{\doi}[1]{\href{https://doi.org/#1}{doi:\nolinkurl{#1}}}
\allowdisplaybreaks
\numberwithin{equation}{section}
\usepackage{subcaption}
\usepackage{float}
\usepackage{multirow}
\usepackage{multicol}
\usepackage{booktabs}
\usepackage{siunitx}

\newtheorem{theorem}{Theorem}[section]
\newtheorem{lemma}[theorem]{Lemma}
\newtheorem{proposition}[theorem]{Proposition}
\newtheorem{corollary}[theorem]{Corollary}
\theoremstyle{definition}

\newtheorem{assumption}[theorem]{Assumption}
\newtheorem{framework}{Framework}

\newtheorem{remark}{Remark}[section]

\newcommand{\bfs}[1]{\boldsymbol{#1}}
\newcommand{\bfx}{\boldsymbol{x}}
\newcommand{\bfv}{\boldsymbol{v}}
\newcommand{\bfy}{\boldsymbol{y}}
\newcommand{\bfz}{\boldsymbol{z}}
\newcommand{\bfn}{\boldsymbol{n}}
\newcommand{\bfu}{\boldsymbol{u}_o}
\newcommand{\bfw}{\boldsymbol{w}}
\newcommand{\bfg}{\boldsymbol{g}}
\newcommand{\bff}{\boldsymbol{f}}
\newcommand{\bfF}{\boldsymbol{F}}

\newcommand{\bfzr}{\boldsymbol{0}}

\renewcommand{\d}{\mathrm{d}}
\newcommand{\A}{\mathcal{A}}
\newcommand{\M}{\mathcal{M}}
\newcommand{\Ueff}{U_{\mathrm{eff}}}

\title{Sea Ice Floe Dynamics and Volume Loss Estimates under Marine Heatwaves}

\newcommand{\FloeAuthorBlock}{%
  \parbox{0.96\textwidth}{\centering
    Quanling Deng\textsuperscript{a,b},\quad
    Seung-Yeal Ha\textsuperscript{c},\quad
    Jaemoon Lee\textsuperscript{d,*},\quad
    Yan Yu Julian Yip\textsuperscript{e}\par
    \vspace{0.7em}
    \small
    \textsuperscript{a}Yau Mathematical Sciences Center, Tsinghua University,
    Jingzhai, Haidian District, Beijing 100084, China\par
    \textsuperscript{b}School of Computing, Australian National University,
    108 North Road, Acton, ACT 2601, Australia\par
    \textsuperscript{c}Department of Mathematical Sciences and Research Institute of Mathematics,
    Seoul National University, 1 Gwanak-ro, Gwanak-gu,
    Seoul 08826, Republic of Korea\par
    \textsuperscript{d}Department of Mathematical Sciences,
    Seoul National University, 1 Gwanak-ro, Gwanak-gu,
    Seoul 08826, Republic of Korea\par
    \textsuperscript{e}Department of Mathematical Sciences, Tsinghua University,
    1 Qinghuayuan, Haidian District, Beijing 100084, China\par
    \vspace{0.5em}
    \textsuperscript{*}Corresponding author:
    \href{mailto:dlwoans0001@snu.ac.kr}{\texttt{dlwoans0001@snu.ac.kr}}%
  }%
}

\author{\FloeAuthorBlock}
\date{}

\begin{document}
\maketitle
\begin{abstract}
We study the ice loss dynamics of an initially observed collection of cylindrical sea ice floes under marine heatwaves. The floe states evolve dynamically under forcings from normal spring--dashpot contact, quadratic ocean and wind drag, Coriolis acceleration, and acceleration due to the slope of the ocean surface. Their radii and thicknesses decrease according to temperature-index melt laws, with basal exchange depending on ocean--ice slip. We prove conditional global existence and uniqueness through extinction under non-rafting and non-grazing contact assumptions. 
For a global solution, we establish upper and lower bounds on the ice loss in two ways. The first method uses environmental data over
admissible regions containing the floe centers. The second evaluates the
melt laws along approximate floe trajectories and bounds the error in
the resulting volume prediction.
Numerical simulations in various scenarios verify the theoretical findings.

\end{abstract}


\section{Introduction}\label{sec:introduction}

Sea ice regulates heat, momentum, freshwater, and biogeochemical exchanges at the air--sea interface
\citep{lund2020arctic,hunke2011multiphase,blockley2020future} in the polar regions. Its mechanical state depends on floe drift,
collisions, and separation under forcings such as Coriolis force, ocean and atmospheric drag forces \citep{parkinson1979large,feltham2008sea,alberello2020drift,alberello2022three}. In the marginal ice zone, the
fragmented cover is mediated by interactions among floes rather than the motion of a coherent sheet
\citep{bennetts2022marginal}. Continuum descriptions, including elastic--plastic, viscous--plastic,
elastic--viscous--plastic, and Maxwell elasto-brittle rheologies, remain important for sea ice dynamics at large scales \citep{coon1974modeling,hibler1979dynamic,hunke1997elastic,dansereau2016maxwell}. In fragmented ice,
however, the dependence on individual contacts and floe size distributions motivates discrete element methods
\citep{hopkins2004discrete,herman2016discrete,damsgaard2018application,manucharyan2022subzero}
and multiscale descriptions linking particle and continuum dynamics
\citep{weiss2017linking,deng2024particle,deng2025particle}.

We study how marine heatwaves change the ice volume of a collection of drifting sea ice floes.
We ask two questions: how the mechanical system evolves through collisions and extinction, and how the ice loss can be bounded during a prescribed marine heatwave without resolving every floe trajectory.
These questions connect the mechanics of fragmented ice with the thermal exposure experienced by its moving constituents.
For simplicity, we adopt a discrete element model (DEM) of non-rotating cylindrical floes with normal contact forces. Such forces
generate no torque with respect to the floe centers, and orientation does not enter their contact geometry. The
resulting spring--dashpot interaction is a reduced discrete element law
\citep{cundall1979discrete,johnson1987contact}. Rotational and tangential effects can be included in more
detailed floe models \citep{damsgaard2018application,deng2026particle}; here they are omitted to isolate the
consequences of lateral recession, basal thinning, and extinction. The intended regime is where size changes due to melting are important and frictional contacts, ridging, and rafting are secondary.
The drag force coefficients are fixed as the floes shrink. The analysis therefore concerns this
specified reduced force law rather than a drag closure that depends on floe area.

Marine heatwaves are prolonged anomalies in ocean temperature relative to a local seasonal climatology. A widely
used definition requires water temperature above a seasonally varying threshold at the 90th percentile for at least five
consecutive days, with gaps of at most two days joined into an event \citep{hobday2016hierarchical}. Their
drivers include oceanic advection, air--sea heat fluxes, and changes in mixing \citep{holbrook2019global}, and
their effects extend to marine ecosystems and coastal economies
\citep{smale2019marine,smith2021socioeconomic}. In the Arctic, heatwave occurrence and persistence have been
associated with declining ice concentration and abrupt sea ice retreat
\citep{huang2021prolonged,barkhordarian2024arctic}. In particular, \cite{huang2021prolonged} documented
increasing intensity, frequency, duration, and areal coverage over 1982--2020, with the season lengthening mainly because it ended later.
These observations motivate a description of how ocean heat melts individual floes during
prescribed warm events.

The temperature driving our melt laws is the local ocean temperature above freezing, not the anomaly above the climatological heatwave threshold. 
We assume that the event and its environmental fields have already been specified. 
Lateral melting changes floe radius and contact geometry, whereas basal melting changes thickness and mass. 
Other thermodynamic sea ice models distinguish between lateral and basal melting \citep{roach2018emergent,bateson2022floe}.
We prescribe melting by ocean heat, omit surface ablation and freezing, and keep the floe interior at a fixed
freezing temperature. Basal heat transfer is represented by an effective exchange speed equal to a positive background contribution plus the relative speed between ocean and ice.
Consequently, mechanics influences melting through both the sampled temperature and basal exchange.

The contact law is discontinuous when the dashpot is active. 
At a non-grazing contact, i.e., a contact at which the gap derivative
is nonzero, the sign of the gap derivative determines the outgoing
contact graph.
Lemma~\ref{lem:wellp} gives local existence and uniqueness while sizes and center distances remain positive. Theorem~\ref{thm:global_continuation} extends the solution through extinction under non-rafting and non-grazing conditions, including limiting contact states. 
Here, non-rafting means that no active center lies strictly inside another floe's disk.
Its proof uses the pointwise velocity bound in Lemma~\ref{lem:pointwise_velocity} to obtain finite limits for the surviving states.

\begin{figure}[t]
    \centering
    \includegraphics[width=\linewidth]{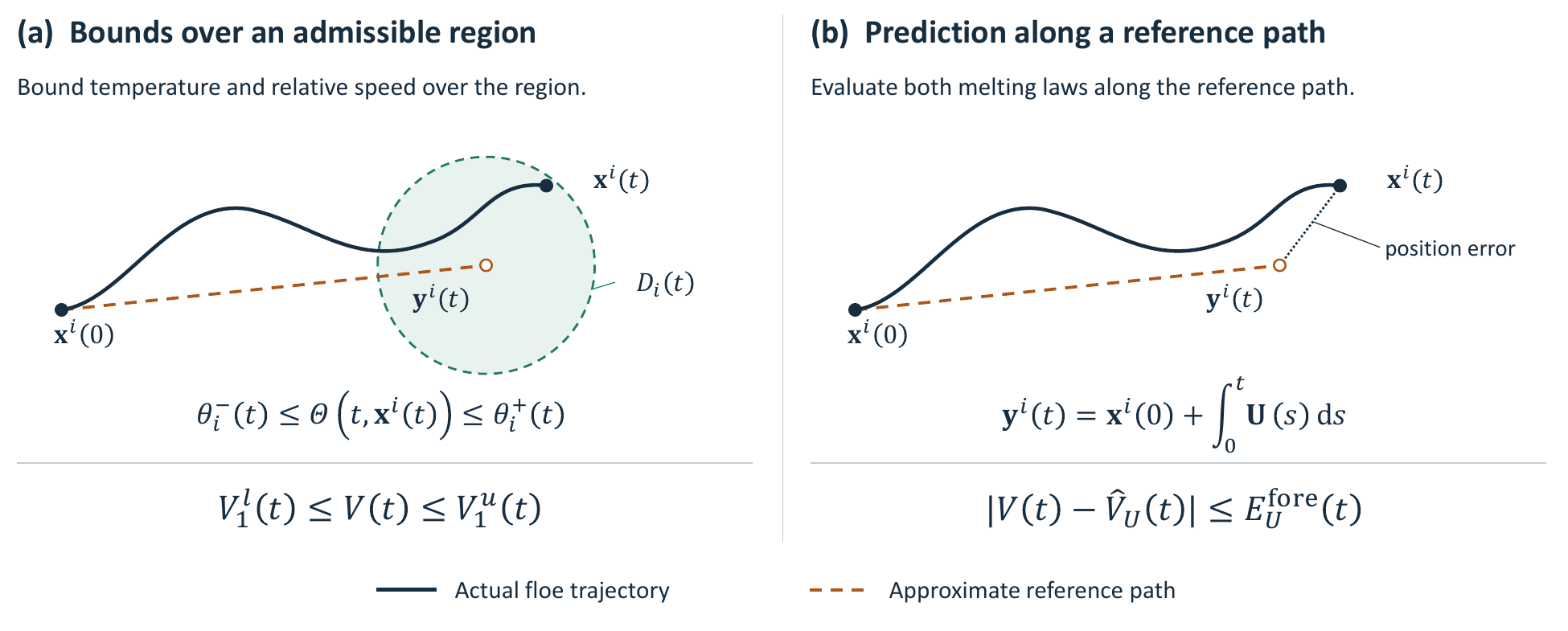}
    \caption{Two approaches to estimating ice loss, illustrated for a single
floe. (a) The first uses bounds on temperature and relative speed over
an admissible region $D_i(t)$ containing the floe center $\bfx^i(t)$.
(b) The second evaluates the melting laws along an approximate trajectory
$\bfy^i(t)$ and bounds the error in the prediction.
The authors used ChatGPT (OpenAI) to help draft this schematic
and verified its scientific content.}
    \label{fig:two_volume_bounds}
\end{figure}

The ice volume estimates use the cubic dissipation supplied by constant quadratic drag. 
Proposition~\ref{prop:motion_budget} bounds the accumulated relative motion using initial states and prescribed fields. Lemma~\ref{lem:weighted_bounds} gives a common H\"older estimate for displacement and melting. 
Theorem~\ref{thm:volume_envelope} then gives lower and upper volume bounds using environmental bounds over admissible regions containing the floe centers, and Theorem~\ref{thm:reference_forecast} bounds the error of a melting prediction along approximate floe trajectories. 
Figure~\ref{fig:two_volume_bounds} illustrates the two constructions
for a single floe.
Corollary~\ref{cor:combined_volume_interval} combines these estimates with an elementary thermal upper bound. 
These estimates concern existing solutions obeying the stated removal rule; no positive temperature or size lower bound is needed.

The rest of the paper is organized as follows. 
Section~\ref{sec:particle} introduces the model. Section~\ref{sec:wellp} establishes local well-posedness and conditional global continuation and defines the aggregate observables.
Section~\ref{sec:prediction} establishes bounds on relative energy and accumulated motion,
develops the ice volume estimates, and discusses their interpretation.
Two numerical experiments are presented in Section \ref{sec:num} to verify the theoretical findings. 
Section~\ref{sec:conclusion} summarizes the conclusions and directions for validation.
The appendices contain the melt calibration and the proofs of the energy estimate and Theorem~\ref{thm:reference_forecast}.

\paragraph*{Notation.}
We write $[n]=\{1,\ldots,n\}$ for the initial labels, with $n\ge1$.
For a scalar $s$, $s_+=\max\{s,0\}$; also
$a\wedge b=\min\{a,b\}$, $\mathbf 1_E$ is the indicator of $E$,
and $\inf\emptyset=\infty$.
A trace $a(\tau^-)$ denotes the left limit of $a$ at $\tau$.

Bold symbols denote planar vectors; $|\cdot|$ and $\cdot$ denote
the Euclidean norm and inner product, respectively.
Matrix norms are the induced operator norms, and
$\|\cdot\|_\infty$ denotes the supremum norm over the stated domain.
For a square-integrable function $f$ on $(0,T)$, we write
\[
 \|f\|_{L^2(0,T)}
 :=\left(\int_0^T |f(t)|^2\,\d t\right)^{1/2}.
\]
The set $\overline B(x,R)$ is the closed ball of radius $R$
centered at $x$.
The rotation $J$ is defined by $J(v_1,v_2)^T=(-v_2,v_1)^T$.

A superscript $i$ denotes a floe label, while $ij$ denotes a pair.
The initial radius and height of floe $i$ are denoted by
$R_i=r^i(0)$ and $H_i=h^i(0)$, respectively.
For a prescribed field $a$, $a^i(t)=a(t,\bfx^i(t))$ denotes
evaluation along the active floe trajectory; the same convention
applies to spatial derivatives.
Unless otherwise specified, label sums at time $t$ are over
$i\in\A(t)$ which denotes the active label set, and pair sums are over $i<j$ with both labels active.

\section{A floe particle model under marine heatwaves}\label{sec:particle}

We specify the coupled mechanical and thermal evolution, its contact and extinction events, and the
observables used in the analysis. The model extends the DEM descriptions of
\cite{damsgaard2018application,herman2016discrete, deng2024particle} for non-rotating floes by allowing each floe to melt both in radius and thickness.

\subsection{Floe states and melting mechanism}\label{sec:floe_states}

We first describe the active state and the two melt rates. Let $\Omega=\mathbb R^2$. Each active floe $i\in[n]$ has state
\[
 (\bfx^i,\bfv^i,r^i,h^i)\in\Omega\times\mathbb R^2\times(0,\infty)^2,
\]
where $\bfx^i$ is the center, $\bfv^i$ is its translational velocity, $r^i$ is its radius, and $h^i$ is its thickness. The mass is determined by the sizes and the constant ice density:
\[
 m^i=\rho_{ice}\pi(r^i)^2h^i.
\]
The environmental forcing is specified by the prescribed fields: the ocean velocity $\bfu$, wind velocity $\bfw$, ocean temperature $T_o$, and acceleration $\bfg$ due to the slope of the ocean surface.

We fix $T_f=-1.9^\circ\mathrm{C}$ as a representative freezing point for seawater with a salinity of $35\,\mathrm{g\,kg^{-1}}$ \citep{shestov2016consolidation} and assume $T_o\ge T_f$.
The nonnegative temperature above freezing and the effective exchange speed are
\[
 \Theta=T_o-T_f,\qquad \Theta_i(t)=\Theta(t,\bfx^i(t)),\qquad \Ueff^i=U_{\min}+|\bfu^i-\bfv^i|,
\]
where $U_{\min}>0$ represents a background velocity scale associated with unresolved turbulent motions at the ice–ocean interface, allowing for enhanced interfacial exchange even under weak resolved relative motion.

Following empirical melting modeling \cite{maykut1987role,steele1992sea,maykut1995solar,nihashi2008bulk} , the lateral and basal melt rates are simplified as 
\begin{equation}\label{eq:meltrh}
 \M_r^i=\varrho_r\Theta_i^{\varsigma_r},\qquad \M_h^i=\varrho_h\Ueff^i\Theta_i^{\varsigma_h}
\end{equation}
with parameters

\begin{equation}\label{eq:rhorh}
 \begin{aligned}
 \varrho_r&=1.6\times10^{-6}\ \mathrm{m\,s^{-1}\,({}^{\circ}C)^{-1.36}},
 &\varsigma_r&=1.36,\\
 \varrho_h&=5.63\times10^{-6}\ \mathrm{({}^{\circ}C)^{-1}},
 &\varsigma_h&=1.
 \end{aligned}
\end{equation}
The lateral melting parameters are estimated in Maykut--Perovich/Steele parameterization \cite{maykut1987role,steele1992sea} and the basal melting parameters are estimated from \cite{hunke1997elastic, maykut1995solar, parkinson1979large, mcphee2003ocean}. 
For general exponents, the units of $\varrho_r$ and $\varrho_h$ are
$\mathrm{m\,s^{-1}\,({}^{\circ}C)^{-\varsigma_r}}$ and $\mathrm{({}^{\circ}C)^{-\varsigma_h}}$, respectively.
Appendix~\ref{app:melt} gives the interpretation in terms of heat flux and the calibration of constants $\varrho_h$ and $U_{\mathrm{min}}$.
The coefficient $\varrho_h$ is independent of the relative speed, whose effect is
incorporated into $\Ueff^i$.

A floe disappears when either radius or thickness reaches zero. 
With size values at this critical time understood as left limits, we define the melt time for floe $i$ as
\[
 \tau^i=\inf\{t>0:r^i(t^-)=0\ \text{or}\ h^i(t^-)=0\},
\]
The dynamical state for floe $i$ is defined on $[0,\tau^i)$. At a finite removal time $\tau=\tau^i$, the radius vanishes first if
$r^i(\tau^-)=0<h^i(\tau^-)$, and the thickness vanishes first if
$h^i(\tau^-)=0<r^i(\tau^-)$; both sizes may also vanish together. All labels attaining zero size at the
same time are removed simultaneously, together with their incident contact forces. When an evolution is considered after removal, its surviving floes are assumed to have finite left limits and to restart from these limits with the same material coefficients. 
We denote $\A(t)$ as the active index set of floes at time $t$. 
This specifies the removal rule for a given evolution; the local existence result below concerns positive sizes only.

\subsection{The floe particle model}

We first define the contact and environmental forces. 
Combining them with the melting laws gives the coupled model. For a pair of floes with distinct centers, we set
\[
 d^{ij}=|\bfx^j-\bfx^i|,\qquad \bfn^{ij}=\frac{\bfx^j-\bfx^i}{d^{ij}},\qquad \delta^{ij}=d^{ij}-(r^i+r^j).
\]
The normal points from floe $i$ to floe $j$, and $\delta^{ij}<0$ means overlap.  
We further define
\[
 q^{ij}=(\bfv^i-\bfv^j)\cdot\bfn^{ij},\qquad \chi^{ij}=\mathbf 1_{\{\delta^{ij}<0\}}.
\]
The normal spring--dashpot force exerted on $i$ by $j$ is
\begin{equation}\label{eq:prmfnft}
 \bff^{ij}=\chi^{ij}(\kappa\delta^{ij}-\gamma q^{ij})\bfn^{ij},\qquad \kappa>0,\quad\gamma\ge0.
\end{equation}
By Newton's third law,$\bff^{ji}=-\bff^{ij}$ and $\bff^{ii}=\bfzr$. The spring is repulsive in overlap, while the signed
dashpot dissipates relative normal motion. During separation the dashpot can produce an attractive resultant; the law is not clipped to enforce unilateral repulsion. 
For $\gamma>0$, since $q^{ij}\bfn^{ij}$ need not vanish at contact onset or offset and $\chi^{ij}$ changes discontinuously, the contact force may exhibit a jump discontinuity at these instants.
So the equations are therefore imposed almost everywhere, and velocities remain continuous at contact switches. 
For $\gamma=0$, the force is locally Lipschitz across contact onset whenever the centers are distinct.

The environmental force consists of quadratic ocean and wind drag, Coriolis acceleration, and acceleration due to the slope of the ocean surface:
\[
 \bfF_e^i=\alpha_o(\bfu^i-\bfv^i)|\bfu^i-\bfv^i|+\alpha_w(\bfw^i-\bfv^i)|\bfw^i-\bfv^i|-m^if_EJ\bfv^i+m^i\bfg^i.
\]
Here $f_E$ is a constant Coriolis parameter and $\bfg(t,\bfx)$ is prescribed
as in \citep{cushman2011introduction,lepparanta2005drift}. 
The drag coefficients $\alpha_o,\alpha_w>0$ 
have units $\mathrm{kg\,m^{-1}}$, and remain constant as the floes evolve. 
For analytical estimates, we denote 
$\alpha:=\alpha_o+\alpha_w>0.$
With the above setting, a coupled system for floe dynamics under a marine heatwave is 
\begin{equation}\label{eq:dem}
\begin{cases}
\displaystyle\dot\bfx^i=\bfv^i,\qquad i\in \A(t),\\
 \displaystyle\dot r^i=-\M_r^i,\\
 \displaystyle\dot h^i=-\M_h^i,\\
 \displaystyle m^i\dot\bfv^i=\sum_{j\in\A(t)\setminus\{i\}}\bff^{ij}+\bfF_e^i,
\end{cases}
\end{equation}
where the dot denotes a time derivative of the quantity. 
We assume that the meltwater leaves with the instantaneous translational velocity of its parent floe. 
Under this assumption, the momentum balance becomes
\[
 \frac{\d}{\d t}(m^i\bfv^i)=\sum_{j\ne i}\bff^{ij}+\bfF_e^i+\dot m^i\bfv^i,
\]
which reduces to the stated equation for $m^i\dot\bfv^i$.

\section{Well-posedness and properties of floe dynamics} \label{sec:wellp}

In this section, we establish local well-posedness and conditional global continuation through floe extinction, and derive balance laws for ice area, volume, and mass.

\subsection{Local existence and uniqueness}

Fix an initial time $t_0\ge0$. We seek a solution forward in time from an initial state at $t_0$. 
On each interval with a fixed active set, a solution of \eqref{eq:dem} is a locally absolutely continuous state function satisfying the equations almost everywhere.
Velocities are continuous at contact switches. At a removal time, surviving floes continue to evolve as \eqref{eq:dem} according to the rule in Section~\ref{sec:floe_states}.
To study the trajectory behavior, we first present the assumptions. 

\begin{assumption}\label{ass:activeAset}
The fields $\bfu,\bfw,\Theta,\bfg$ belong to $C_b^1([0,T]\times\Omega)$, with their respective scalar or vector values, for every $T>0$, and $\Theta\ge0$.  
The initial sizes of the floes are positive, and their initial centers are pairwise distinct. 
\end{assumption}
The melt coefficients are positive and $\varsigma_r,\varsigma_h\ge1$, as in \eqref{eq:rhorh}. These exponent restrictions make the temperature powers locally Lipschitz at zero and ensure the regularity used below.
When two floes collide, taking the time derivative of the overlap $\delta^{ij}$ leads to,
\begin{equation}\label{eq:gap_derivative}
 \dot\delta^{ij}=-q^{ij}+\M_r^i+\M_r^j.
\end{equation}
With this in mind, we present the non-grazing assumption below.

\begin{assumption}\label{ass:regularcontacts}
If $\gamma>0$, every floe pair in contact at $t_0$ is non-grazing, that is,
\[
 \delta^{ij}(t_0)=0
 \quad\Longrightarrow\quad
 -q^{ij}(t_0)+\M_r^i(t_0)+\M_r^j(t_0)\ne0.
\]
\end{assumption}

The next lemma establishes local existence and uniqueness, including simultaneous contact switches.
\begin{lemma}\label{lem:wellp}
Suppose Assumption~\ref{ass:activeAset} holds. Every state with positive radii and thicknesses
and distinct centers satisfying Assumption~\ref{ass:regularcontacts}
at $t_0$ admits a unique solution of \eqref{eq:dem} on
$[t_0,t_0+\varepsilon)$ for some $\varepsilon>0$.
\end{lemma}
\begin{proof}
First we solve the equations with a fixed contact graph. Second we select the outgoing branch.

\noindent\textit{$\bullet$ Step 1} (Fixed contact graph).
Choose a candidate contact graph containing all pairs with $\delta^{ij}(t_0)<0$ and a subset of the pairs with $\delta^{ij}(t_0)=0$.
In the auxiliary equations,
replace $\chi^{ij}$ by $1$ on the edges of this graph and by $0$ on all other
pairs, independently of the evolving gaps. With these choices held fixed,
positive masses and distinct centers make the right-hand side continuous in
time and locally Lipschitz in the state. Picard--Lindel\"of gives a unique
branch solution. Step 2 checks that its gaps agree with the chosen graph.
If $\gamma=0$, the force
$\bff^{ij}=\kappa\min\{\delta^{ij},0\}\bfn^{ij}$ is locally Lipschitz
across contact as well, proving the assertion directly.

\noindent\textit{$\bullet$ Step 2} (Contact switches).
Let $\gamma>0$. Retain pairs with negative gap and include each colliding pair
exactly when $\dot\delta^{ij}(t_0)<0$. By \eqref{eq:gap_derivative}, all
branches have the same nonzero derivative at a colliding pair. The chosen
branch therefore satisfies
\[
 \delta^{ij}(t)=\dot\delta^{ij}(t_0)(t-t_0)+o(t-t_0),
 \qquad t\downarrow t_0.
\]
Thus all chosen signs hold on a common short interval, and initially nonzero
gaps keep their signs. The indicator value at $t_0$ does not affect the
integral equation. Every solution has the same initial gap derivatives by
the position and radius equations, so it follows this branch. Step 1 proves
uniqueness.
\end{proof}

\subsection{Conditional global continuation through extinction}
\label{subsec:global_continuation}
We first obtain a pointwise velocity bound from the quadratic drag.
We then use this bound to prove conditional continuation through extinction in Theorem~\ref{thm:global_continuation}.

Fix $T<\infty$. For $0\le t\le T$, let
\[
 v_{\max}(t):=\max_{i\in\A(t)}|\bfv^i(t)|,
 \qquad R_*:=\max_i R_i,\qquad \beta:=2(n-1)\gamma,
\]
with $v_{\max}(t)=0$ if $\A(t)=\emptyset$. 
Let $V_{\mathrm{crit}}$ be the largest real root of the following quadratic equation in $V$:
\[
 \alpha\bigl(V-\max\{\|\bfu\|_\infty,\|\bfw\|_\infty\}\bigr)^2
 =\beta V+2(n-1)\kappa R_*+\|\bfg\|_\infty\max_i m^i(0),
\]
where the field norms are defined over $[0,T]\times\Omega$.

\begin{lemma}[Pointwise velocity bound]\label{lem:pointwise_velocity}
Suppose that a solution satisfying the removal rule is defined on an interval $I\subset[0,T]$ with left endpoint $s\in I$.
Then every floe active at $t\in I$ satisfies
\[
 |\bfv^i(t)|\le\max\{v_{\max}(s),V_{\mathrm{crit}}\}.
\]
\end{lemma}
\begin{proof}
Suppose $v_{\max}(t)>V_{\mathrm{crit}}$ and let $i$ be a label attaining
$v_{\max}(t)$. Since
$V_{\mathrm{crit}}\ge\max\{\|\bfu\|_\infty,\|\bfw\|_\infty\}$,
for each $\bfs a\in\{\bfu^i,\bfw^i\}$ we have
\[
 (\bfv^i-\bfs a)\cdot\frac{\bfv^i}{v_{\max}}
 \ge v_{\max}-|\bfs a|,\qquad
 |\bfv^i-\bfs a|\ge v_{\max}-|\bfs a|.
\]
Each contact force has magnitude at most $2(\kappa R_*+\gamma v_{\max})$
by \eqref{eq:prmfnft}, and the Coriolis force is orthogonal to $\bfv^i$.
Since index $i$ is a maximizing label, $\dot v_{\max}=\frac{\d}{\d t}|\bfv^i|$ almost everywhere.
Thus it follows from \eqref{eq:dem} that
\[
 \begin{aligned}
 m^i\frac{\d}{\d t}v_{\max}
 &\le -\alpha\bigl(v_{\max}
          -\max\{\|\bfu\|_\infty,\|\bfw\|_\infty\}\bigr)^2\\
 &\quad +\beta v_{\max}+2(n-1)\kappa R_*
          +\bigl(\max_j m^j(0)\bigr)\|\bfg\|_\infty
 <0,
 \end{aligned}
\]
almost everywhere on $\{v_{\max}>V_{\mathrm{crit}}\}$.
Since $v_{\max}$ is locally absolutely continuous and $m^i>0$ on the
phase, the bound follows. 
At each floe removal, the velocities of the remaining floes are continuous across the removal time, so the same bound applies after restart.
\end{proof}

Recall that a configuration is non-rafting if no active floe center lies strictly inside the disk of another floe. The following framework formalizes this geometric condition together with the additional assumptions required for conditional continuation.

\begin{framework}[Conditional global continuation]\label{fw:global_continuation}
The locally constructed solution is assumed to satisfy the following conditions throughout its interval of existence:
\begin{enumerate}
\item[(i)] Assumption~\ref{ass:activeAset} holds.
\item[(ii)] Along each phase of the locally constructed evolution with positive radii and thicknesses,
\begin{equation}\label{eq:non-rafting}
 d^{ij}(t)\ge\max\{r^i(t),r^j(t)\},\qquad i,j\in\A(t),\quad i\ne j.
\end{equation}
\item[(iii)] If $\gamma>0$, Assumption~\ref{ass:regularcontacts} holds whenever two active floes collide, that is, when $\delta^{ij}=0$, and at every finite limiting contact state of labels whose limiting radii and thicknesses are positive.
\end{enumerate}
\end{framework}
Conditions (ii)--(iii) are imposed a priori along the evolution, not only at the initial time;
their preservation is not proved here.
Thus Theorem~\ref{thm:global_continuation} is a conditional continuation result.
For $\gamma=0$, condition (iii) is unnecessary.

\begin{theorem}[Global continuation through extinction]\label{thm:global_continuation}
Under Framework~\ref{fw:global_continuation}, the following assertions hold:
\begin{enumerate}
\item (Global existence and uniqueness)
The local solution of \eqref{eq:dem} admits a unique global continuation under the removal rule
in Section~\ref{sec:floe_states}.
\item (Behavior at extinction)
There are at most $n$ removal times. At each removal time, every surviving floe continues from the
finite limit of its full state, while every disappearing floe has a finite position limit.
\end{enumerate}
Extinction may occur through the vanishing of either the radius or the thickness, or both simultaneously.
\end{theorem}
\begin{proof}
We first establish the endpoint limits in assertion (2), then construct the continuation in assertion (1)
and bound the number of removals. Fix $T<\infty$ and consider a phase with fixed active set $I$,
beginning at $s$ and ending at a finite maximal time $t_*\le T$.

\noindent\textit{$\bullet$ Step 1} (Endpoint limits for assertion (2)).
By Lemma~\ref{lem:pointwise_velocity},
$|\bfv^i(t)|\le V_*:=\max\{v_{\max}(s),V_{\mathrm{crit}}\}$
throughout the phase. Hence
\[
 |\bfx^i(t)-\bfx^i(t')|\le V_*|t-t'|,
 \qquad s\le t',t<t_*,
\]
so every center has a finite limit at $t_*$. The nonincreasing positive
sizes also have limits $r_*^i,h_*^i\ge0$.

Let $I_+=\{j\in I:r_*^jh_*^j>0\}$. For $j\in I_+$, the mass is bounded
below near $t_*$. The contact law \eqref{eq:prmfnft} gives
$|\bff^{ji}(t)|\le 2(\kappa R_*+\gamma V_*)$.
The environmental force is also bounded, so \eqref{eq:dem} gives
$|\dot{\bfv}^j(t)|\le C_j$ almost everywhere near $t_*$.
Hence $\bfv^j$ has a finite limit at $t_*$.
Every retained label therefore has a finite limit of its
full state, and \eqref{eq:non-rafting} keeps retained centers distinct. 
This proves the claims about endpoint limits in assertion (2).

\noindent\textit{$\bullet$ Step 2} (Global continuation and the number of removals).
If $I_+=I$, all limiting sizes are positive. For $\gamma>0$, condition (iii) of
Framework~\ref{fw:global_continuation} gives, at each limiting zero gap,
\[
 \dot\delta^{ij}(t_*^-)
 =-q^{ij}(t_*^-)+\M_r^i(t_*^-)+\M_r^j(t_*^-)\ne0.
\]
Its derivative therefore has a fixed sign sufficiently near $t_*$, so the gap is strictly monotone
there and contact switches cannot accumulate at that endpoint. Lemma~\ref{lem:wellp} restarts the
solution from its limiting state, contradicting maximality of the phase. For $\gamma=0$, the
locally Lipschitz elastic field gives the same continuation directly.

Thus a genuine finite phase endpoint has $I\setminus I_+\ne\varnothing$. Remove these labels
simultaneously. The limiting state on $I_+$ has positive sizes, distinct centers, and, when needed,
the non-grazing contact signs, so Lemma~\ref{lem:wellp} gives a unique restart.
The retained velocities at restart still satisfy $|\bfv^j(t_*)|\le V_*$,
so Lemma~\ref{lem:pointwise_velocity} gives the same bound on the next phase.
Each removal deletes at least one label permanently, so there are at most $n$ removal times.
Repeating the construction on every finite horizon proves assertion (1) and completes the proof of assertion (2).
If no active labels remain, all floes have melted.
\end{proof}

\subsection{Balance laws for area, volume, and mass}\label{sec:observables}

We define the aggregate observables and derive the balances connecting mechanical dissipation with the evolution of volume. 
We use physical area and volume throughout:
\begin{equation}\label{eq:IceC}
 A(t)=\pi\sum_{i\in\A(t)}(r^i)^2,\qquad
 V(t)=\pi\sum_{i\in\A(t)}(r^i)^2h^i.
\end{equation}
These quantities track the initially observed floes under the assumption that no floes enter or leave the observation window. Since the floes may overlap, $A$ is the sum of their cross-sectional areas, not the area of their union. 
If desired, $A/A_{\mathrm{ref}}$ and $V/(A_{\mathrm{ref}}h_{\mathrm{ref}})$ give normalized observables for fixed positive reference area and thickness; the estimates below are divided by the same constants.

Between removals, the melt laws give
\[
 \dot m^i=-\rho_{ice}\pi\bigl(2r^ih^i\M_r^i+(r^i)^2\M_h^i\bigr)\le0,
\]
\[
 \dot A=-2\pi\sum_i r^i\M_r^i,\qquad
 \dot V=-\pi\sum_i\bigl(2r^ih^i\M_r^i+(r^i)^2\M_h^i\bigr).
\]
At a simultaneous removal time $\tau$,
\begin{equation}\label{eq:concentration_jumps}
 A(\tau)-A(\tau^-)=-\pi\sum_{i:\tau^i=\tau}r^i(\tau^-)^2,
 \qquad V(\tau)=V(\tau^-).
\end{equation}
Thus, the floe extinction due to radius melting preserves the continuity of sea ice concentration, whereas the extinction due to thickness may induce a downward jump.
Both observables are nonincreasing. 
Integrating the melting laws over successive intervals between floe removals and using the continuity of volume at removal times, we conclude that \(V\) is absolutely continuous on every finite time interval.
In contrast, the floe total area $A$ may exhibit jumps at removal times,
as described in \eqref{eq:concentration_jumps}.

\section{Estimates for the ice loss}\label{sec:prediction}

In this section, we derive estimates for the ice loss
in terms of the initial data and prescribed environmental fields,
without computing the actual floe trajectories.
We first derive an elementary thermal upper bound from the melt laws.
In Theorem~\ref{thm:volume_envelope}, we then derive upper and lower
bounds on the volume using temperature bounds over regions containing
the floe centers.
Finally, in Theorem~\ref{thm:reference_forecast}, we estimate the error
in the volume obtained by evaluating the melt laws along approximate floe trajectories.

Throughout this section we consider a solution on $[0,T]$ satisfying the removal rule. 
Theorem~\ref{thm:global_continuation} supplies such a solution
under Framework~\ref{fw:global_continuation}.
Its volume is defined in \eqref{eq:IceC} and is continuous at each removal by \eqref{eq:concentration_jumps}.

Integrating the two melt laws up to the time of removal gives the volume directly:
\begin{equation}\label{eq:two_exposure_outputs}
 V(t)=\pi\sum_{i=1}^n
 \left(R_i-\int_0^{t\wedge\tau^i}\M_r^i(s)\,\d s\right)_+^2
 \left(H_i-\int_0^{t\wedge\tau^i}\M_h^i(s)\,\d s\right)_+.
\end{equation}
Before removal, the two positive parts are the actual radius and thickness. At removal, at least
one factor vanishes, so the product remains zero thereafter.

An elementary thermal upper bound follows immediately. 
Suppose that $\theta_i^-(s)$ are bounded, measurable, and nonnegative functions satisfying $\theta_i^-(s)\le\Theta(s,\bfx^i(s))$ whenever label $i$ is active.
For instance, one may use a spatially uniform nonnegative lower bound on the temperature above freezing, $\Theta=T_o-T_f$, for all $i$.
Then
\begin{equation}\label{eq:thermal_volume_upper}
 \begin{aligned}
 V(t)\le V_{\mathrm{th}}(t):={}&\pi\sum_{i=1}^n
 \left(R_i-\varrho_r\int_0^t\theta_i^-(s)^{\varsigma_r}\,\d s\right)_+^2
 \left(H_i-\varrho_h U_{\min}\int_0^t\theta_i^-(s)^{\varsigma_h}\,\d s\right)_+.
 \end{aligned}
\end{equation}

Indeed, \eqref{eq:meltrh} and $\Ueff^i\ge U_{\min}$ give the corresponding lower bounds for both melt rates while label $i$ is active. Integrating them bounds its remaining volume from above; after removal its contribution is zero and the comparison remains valid.

\begin{remark}[Pointwise velocity bounds]\label{rem:pointwise_velocity}
To derive a corresponding lower bound on the ice volume by the same
pointwise comparison, we would need an upper bound on the relative speed
between each floe and the ocean.
The pointwise velocity bound of Lemma~\ref{lem:pointwise_velocity}
provides such a bound via the triangle inequality, but using it in the
basal melt integral may yield an overly conservative lower bound on the
remaining ice volume.
For example, the parameters of Example~1 in Section~\ref{sec:num} give
$V_{\mathrm{crit}}\ge\beta/\alpha\simeq 3.6\times10^8\,\mathrm{m\,s^{-1}}$.
Consequently, the resulting lower bound on the ice volume vanishes within
a very short time.
We therefore use an integral bound on relative motion to control accumulated basal melting in Theorem~\ref{thm:volume_envelope}, rather than inserting a pointwise velocity bound.
\end{remark}

\subsection{Relative energy and integrated relative speed}

The main goal of this subsection is to control the deviations of the
floe positions and velocities from a prescribed reference motion.
To this end, we choose a reference velocity and estimate the associated
force residual in terms of the initial data and environmental fields.
More precisely, we aim to derive upper bounds on
\[
 \int_0^t\sum_{i\in\A(s)}|\bfv^i(s)-\bfs U(s)|^3\,\d s,\qquad
 \int_0^t\sum_{i\in\A(s)}K_i(s)|\bfv^i(s)-\bfs U(s)|\,\d s,
\]
for any nonnegative $K_i\in L^{3/2}(0,t)$.

Choose a prescribed, spatially uniform reference velocity $\bfs U\in C^1([0,T];\mathbb R^2)$.
The choice of $\bfs U$ is arbitrary.  
For instance, one may take the zero velocity, a representative wind or ocean velocity, or any other prescribed reference velocity satisfying the above regularity.
We define the relative energy with respect to $\bfs U$ as follows:
\begin{equation}\label{eq:relative_energy}
 E_U(t)=\frac12\sum_{i\in\A(t)}m^i|\bfv^i-\bfs U|^2
          +\frac\kappa2\sum_{i<j}(-\delta^{ij})_+^2.
\end{equation}
Recall $\alpha=\alpha_o+\alpha_w$ and set $\Phi(\bfz):=|\bfz|\bfz$. 
We define the reference force residual by evaluating the environmental
force at the reference velocity and subtracting the reference acceleration
term:
\begin{equation}\label{eq:reference_force_residual}
 \begin{split}
 \bfs B_U(t,\bfx,m)={}&\alpha_o \Phi(\bfu(t,\bfx)-\bfs U(t))
  +\alpha_w \Phi(\bfw(t,\bfx)-\bfs U(t))\\
 &+m\bigl(\bfg(t,\bfx)-\dot{\bfs U}(t)-f_EJ\bfs U(t)\bigr).
 \end{split}
\end{equation}
Note that $\bfs B_U(t,\bfx,m)$ is affine in $m$.
Since $0\le m^i(t)\le m^i(0)$, the vector
$\bfs B_U(t,\bfx^i(t),m^i(t))$ lies on the line segment joining
its values at $m=0$ and $m=m^i(0)$.
Hence, by the convexity of the norm, we obtain
\[
 |\bfs B_U(t,\bfx^i(t),m^i(t))|
 \le \sup_{\bfx\in\Omega}\max\bigl\{
       |\bfs B_U(t,\bfx,0)|,
       |\bfs B_U(t,\bfx,m^i(0))|\bigr\}=: b_i(t).
\]
In particular, $b_i(t)$ can be computed from the initial data,
the environmental fields, and the prescribed reference velocity
without knowledge of the floe trajectories.
We then define the following quantity independent of the trajectory:
\begin{equation}\label{eq:explicit_relative_motion_budget}
 R_U(t)=\frac{3E_U(0)}{\alpha}
       +\frac{2\sqrt2}{\alpha^{3/2}}
          \int_0^t\sum_{i=1}^n b_i(s)^{3/2}\,\d s.
\end{equation}

We are now ready to derive an upper bound on the time integral
of the sum of cubed velocity deviations.

\begin{proposition}\label{prop:motion_budget}
Under the above assumptions and setting, the solution of \eqref{eq:dem} satisfies, for $0\le t\le T$,
\begin{equation}\label{eq:weighted_speed_result}
 \frac{3}{\alpha}E_U(t)
 +\int_0^t\sum_{i\in\A(s)}
   \left(
     |\bfv^i(s)-\bfs U(s)|^3
     +\frac{3}{2\alpha}(-\dot m^i(s))|\bfv^i(s)-\bfs U(s)|^2
   \right)\,\d s
 \le R_U(t).
\end{equation}
\end{proposition}
\begin{proof}
Since the proof is lengthy, we defer it to Appendix~\ref{app:motion_budget_explicit}.
In the proof, we first derive the relative energy identity, then apply Young's inequality to the residual force term, and finally integrate the resulting inequality while accounting for floe removal.
\end{proof}
In particular, since $E_U(t)\ge0$ and $\dot m^i(s)\le0$, we obtain
\begin{equation}\label{eq:slip_budget_for_motion}
 \int_0^t\sum_{i\in\A(s)}|\bfv^i(s)-\bfs U(s)|^3\,\d s\le R_U(t).
\end{equation}

We next apply H\"older's inequality to obtain a weighted estimate
from the bound on the time integral of cubed velocity deviations.
We use this estimate to control the displacement of individual floes,
derive volume bounds in Theorem~\ref{thm:volume_envelope},
and estimate the error of the melting prediction along approximate floe trajectories in Theorem~\ref{thm:reference_forecast}.

\begin{lemma}\label{lem:weighted_bounds}
For any nonnegative weights $K_i\in L^{3/2}(0,t)$, $i\in [n]$,
we have
\begin{equation}\label{eq:weighted_motion_bound}
 \int_0^t\sum_{i\in\A(s)}K_i(s)|\bfv^i(s)-\bfs U(s)|\,\d s
 \le R_U(t)^{1/3}
       \left(\sum_{i=1}^n\int_0^t K_i(s)^{3/2}\,\d s\right)^{2/3}.
\end{equation}
\end{lemma}
\begin{proof}
H\"older's inequality on time and labels gives
\[
 \begin{aligned}
 \int_0^t\sum_{i\in\A(s)}K_i|\bfv^i-\bfs U|\,\d s
 &\le\left(\int_0^t\sum_{i\in\A(s)}|\bfv^i-\bfs U|^3\,\d s\right)^{1/3}
       \left(\int_0^t\sum_{i\in\A(s)}K_i^{3/2}\,\d s\right)^{2/3}\\
 &\le R_U(t)^{1/3}
       \left(\sum_{i=1}^n\int_0^t K_i^{3/2}\,\d s\right)^{2/3}.
 \end{aligned}
\]
\end{proof}
Taking $K_i\equiv1$ for a fixed label $i$ and $K_j\equiv0$ for $j\ne i$ gives, for $0\le t\le T$ with $t<\tau^i$,
\begin{equation}\label{eq:displacement_bound}
 \left|\bfx^i(t)-\bfx^i(0)-\int_0^t\bfs U(s)\,\d s\right|
 \le\int_0^t|\bfv^i(s)-\bfs U(s)|\,\d s
 \le t^{2/3}R_U(t)^{1/3}.
\end{equation}
Thus, for $t<\tau^i$, the floe position $\bfx^i(t)$ lies in the
closed ball centered at $\bfx^i(0)+\int_0^t\bfs U(s)\,\d s$
with radius $t^{2/3}R_U(t)^{1/3}$.

\subsection{Melting prediction over admissible regions}

Here we present the first of two estimates that account for the effect
of floe motion on the ice loss.
The strategy is as follows.
First, for each active floe, \eqref{eq:displacement_bound} gives a region
containing its center.
Then, we bound the temperature along its trajectory by taking the infimum and
the supremum of the prescribed lower/upper temperature bound over this region.
Finally, we compare the ocean--ice relative speed $|\bfu^i-\bfv^i|$
with $|\bfu^i-\bfs U|$. Lemma~\ref{lem:weighted_bounds} controls the
accumulated contribution of the velocity difference $|\bfv^i-\bfs U|$.

We start by introducing terms that are needed to state the volume bound.
Suppose the uncertain temperature above freezing satisfies prescribed bounded continuous bounds
\[
 0\le\underline\Theta(t,\bfx)\le\Theta(t,\bfx)\le\overline\Theta(t,\bfx).
\]
For each $i\in[n]$, choose nonempty prescribed regions $D_i(t)\subset\Omega$ for every $t\in[0,T]$, with $\bfx^i(t)\in D_i(t)$ whenever label $i$ is active. 
One possible choice is
\[
 D_i(t)=\overline B\left(\bfx^i(0)+\int_0^t\bfs U(s)\,\d s,\
                        t^{2/3}R_U(t)^{1/3}\right).
\]
The inclusion $\bfx^i(t)\in D_i(t)$ for $t<\tau^i$ follows from
\eqref{eq:displacement_bound}.
The global choice $D_i(t)=\Omega$ is also admissible.
The regions must be chosen so that the infima and suprema defining $\theta_i^\pm$ and $c_i^\pm$ below are measurable functions of time. 
By Berge's maximum theorem~\cite{berge1963topological}, this condition is automatically satisfied for compact regions varying continuously in the Hausdorff metric, since the
underlying fields are continuous.
In particular, it holds for the closed balls specified above.

From these regions and the prescribed fields, compute
\begin{equation}\label{eq:tube_temperatures}
 \theta_i^-(t)=\inf_{\bfx\in D_i(t)}\underline\Theta(t,\bfx),\qquad
 \theta_i^+(t)=\sup_{\bfx\in D_i(t)}\overline\Theta(t,\bfx).
\end{equation}
Define the radius bounds
\[
 \mathfrak l_i(t)=\left(R_i-\varrho_r\int_0^t\theta_i^+(s)^{\varsigma_r}\,\d s\right)_+,
 \qquad
 \mathfrak u_i(t)=\left(R_i-\varrho_r\int_0^t\theta_i^-(s)^{\varsigma_r}\,\d s\right)_+.
\]
They provide lower and upper bounds on the radius of floe $i$:
$\mathfrak l_i(t)\le r^i(t)\le\mathfrak u_i(t)$
whenever label $i$ is active.

For each $i\in[n]$, define
\[
 w_i(s):=1+\frac{3\rho_{ice}\pi\varrho_h}{2\alpha}
                 \mathfrak l_i(s)^2\theta_i^-(s)^{\varsigma_h}.
\]
We choose the weight $w_i$ to make use of the integrand in
\eqref{eq:weighted_speed_result} of Proposition~\ref{prop:motion_budget};
see \eqref{eq:weighted_excess_integrand}.

We bound the relative speed between the ocean and the reference velocity over $D_i(s)$ and add $U_{\min}$ to obtain 
\[
 c_i^-(s)=U_{\min}+\inf_{\bfx\in D_i(s)}|\bfu(s,\bfx)-\bfs U(s)|,\qquad
 c_i^+(s)=U_{\min}+\sup_{\bfx\in D_i(s)}|\bfu(s,\bfx)-\bfs U(s)|.
\]

All environmental integrals below extend over $[0,T]$, even when a label disappears earlier. 
If the ocean, wind, or slope fields are uncertain, the bounds used to construct the volume interval must hold uniformly over all admissible fields.

\begin{theorem}[Volume bounds from reachable regions]\label{thm:volume_envelope}
Under the preceding assumptions, define
\begin{equation}\label{eq:forecast_volume_bounds}
 \begin{aligned}
 V_1^l(t)=\pi\Biggl[&\sum_{i=1}^n\mathfrak l_i(t)^2
      \left(H_i-\varrho_h\int_0^t
          c_i^+(s)\theta_i^+(s)^{\varsigma_h}\,\d s\right)_+-\varrho_h R_U(t)^{1/3}
      \left(\sum_{i=1}^n\mathfrak l_i(t)^3
          \int_0^t\frac{\theta_i^+(s)^{3\varsigma_h/2}}
                        {\sqrt{w_i(s)}}\,\d s\right)^{2/3}\Biggr]_+,\\
 V_1^u(t)=\pi\Biggl[&\sum_{i=1}^n\mathfrak u_i(t)^2
      \left(H_i-\varrho_h\int_0^t c_i^-(s)\theta_i^-(s)^{\varsigma_h}\,\d s\right)_+ +\varrho_h R_U(t)^{1/3}
      \left(\sum_{i=1}^n\mathfrak u_i(t)^3
              \int_0^t\theta_i^-(s)^{3\varsigma_h/2}\,\d s\right)^{2/3}\Biggr].
 \end{aligned}
\end{equation}
For $0\le t\le T$, 
\begin{equation}\label{eq:forecast_volume_envelope}
 V_1^l(t)\le V(t)\le V_1^u(t).
\end{equation}

\end{theorem}
\begin{proof}
We first compare each floe's volume and then sum the resulting velocity errors.

\medskip
\noindent\textbf{$\bullet$ Step 1} (Individual volume bounds).
While label $i$ is active, define the excess above $c_i^+$ and the
shortfall below $c_i^-$ by
\[
 p_i(s):=\bigl(\Ueff^i(s)-c_i^+(s)\bigr)_+,\qquad
 q_i(s):=\bigl(c_i^-(s)-\Ueff^i(s)\bigr)_+.
\]
Since $c_i^-\le U_{\min}+|\bfu^i-\bfs U|\le c_i^+$, the reverse
triangle inequality gives
\[
 \max\{p_i,q_i\}
 \le \bigl||\bfu^i-\bfv^i|-|\bfu^i-\bfs U|\bigr|
 \le |\bfv^i-\bfs U|.
\]
Moreover, $c_i^+\ge U_{\min}$ gives
$p_i=(U_{\min}+|\bfu^i-\bfv^i|-c_i^+)_+\le|\bfu^i-\bfv^i|$.
Thus
\begin{equation}\label{eq:exchange_error_bounds}
 p_i\le\min\{|\bfv^i-\bfs U|,\,|\bfu^i-\bfv^i|\},
 \qquad
 q_i\le|\bfv^i-\bfs U|.
\end{equation}
By definition,
\[
 0\le c_i^-(s)-q_i(s)\le \Ueff^i(s)\le c_i^+(s)+p_i(s).
\]
Combining these bounds with the temperature bounds in the basal melt law
yields
\[
 \varrho_h\theta_i^-(s)^{\varsigma_h}\bigl(c_i^-(s)-q_i(s)\bigr)
 \le \M_h^i(s)
 \le \varrho_h\theta_i^+(s)^{\varsigma_h}\bigl(c_i^+(s)+p_i(s)\bigr).
\]
For fixed $t$, define the accumulated errors
\[
 \begin{aligned}
 e_i^+&:=\varrho_h\int_0^{t\wedge\tau^i}
          \theta_i^+(s)^{\varsigma_h}p_i(s)\,\d s,\\
 e_i^-&:=\varrho_h\int_0^{t\wedge\tau^i}
          \theta_i^-(s)^{\varsigma_h}q_i(s)\,\d s.
 \end{aligned}
\]
For $t<\tau^i$, integration of the basal melt bounds, together with
$\mathfrak l_i(t)\le r^i(t)\le\mathfrak u_i(t)$, gives
\[
 \begin{aligned}
 (r^i(t))^2h^i(t)&\ge \mathfrak l_i(t)^2
   \left(H_i-\varrho_h\int_0^t
     c_i^+(s)\theta_i^+(s)^{\varsigma_h}\,\d s-e_i^+\right)_+,\\
 (r^i(t))^2h^i(t)&\le \mathfrak u_i(t)^2
   \left(H_i-\varrho_h\int_0^t
     c_i^-(s)\theta_i^-(s)^{\varsigma_h}\,\d s+e_i^-\right)_+.
 \end{aligned}
\]
These inequalities also hold after removal, with the actual contribution set to zero. 
The upper bound remains valid by nonnegativity.
For the lower bound, if radius vanishes at $\tau^i$, the radial comparison gives $\mathfrak l_i(\tau^i)=0$, and this remains zero. If thickness vanishes instead, the integrated upper bound on basal melting gives
\[
 H_i-\varrho_h\int_0^{\tau^i}c_i^+(s)\theta_i^+(s)^{\varsigma_h}\,\d s-e_i^+\le0.
\]
For later $t$, $e_i^+$ is frozen at removal and the nonnegative baseline integral can only increase, so the lower product remains zero. 

Finally, $(a-e)_+\ge a_+-e$ and $(a+e)_+\le a_++e$ for $e\ge0$
allow us to take $e_i^+$ and $e_i^-$ outside the positive parts.

\medskip
\noindent\textbf{$\bullet$ Step 2} (Summation using the bound on accumulated relative motion).
While label $i$ is active, dropping the nonnegative lateral melt term
and using \eqref{eq:exchange_error_bounds} gives
\begin{equation}\label{eq:mass_loss_excess}
 \begin{aligned}
 -\dot m^i
 &=\rho_{ice}\pi\bigl(2r^ih^i\M_r^i+(r^i)^2\M_h^i\bigr)\\
 &\ge \rho_{ice}\pi\varrho_h\mathfrak l_i^2
       (\theta_i^-)^{\varsigma_h}\Ueff^i\\
 &\ge \rho_{ice}\pi\varrho_h\mathfrak l_i^2
       (\theta_i^-)^{\varsigma_h}p_i
  =\frac{2\alpha}{3}(w_i-1)p_i.
 \end{aligned}
\end{equation}
Multiplying \eqref{eq:mass_loss_excess} by $3p_i^2/(2\alpha)$ and
using $p_i\le|\bfv^i-\bfs U|$, we obtain
\begin{equation}\label{eq:weighted_excess_integrand}
 \begin{aligned}
 w_i p_i^3
 &\le p_i^3+\frac{3}{2\alpha}(-\dot m^i)p_i^2\\
 &\le |\bfv^i-\bfs U|^3
      +\frac{3}{2\alpha}(-\dot m^i)|\bfv^i-\bfs U|^2.
 \end{aligned}
\end{equation}
Hence dropping nonnegative $E_U(t)$ from Proposition~\ref{prop:motion_budget} yields
\[
 \int_0^t\sum_{i\in\A(s)}w_i(s)p_i(s)^3\,\d s\le R_U(t).
\]
For fixed $t$, applying H\"older's inequality to
$w_i^{1/3}p_i$ and
$\mathfrak l_i(t)^2\theta_i^+(s)^{\varsigma_h}w_i^{-1/3}$ gives
\[
 \begin{aligned}
 \sum_{i=1}^n\mathfrak l_i(t)^2e_i^+
 &=\varrho_h\int_0^t\sum_{i\in\A(s)}
      \mathfrak l_i(t)^2\theta_i^+(s)^{\varsigma_h}p_i(s)\,\d s\\
 &\le \varrho_h R_U(t)^{1/3}
   \left(\sum_{i=1}^n\mathfrak l_i(t)^3
      \int_0^t\frac{\theta_i^+(s)^{3\varsigma_h/2}}
                    {\sqrt{w_i(s)}}\,\d s\right)^{2/3}.
 \end{aligned}
\]
For the upper bound, Lemma~\ref{lem:weighted_bounds} with
$K_i(s)=\varrho_h\mathfrak u_i(t)^2\theta_i^-(s)^{\varsigma_h}$ gives
\[
 \begin{aligned}
 \sum_{i=1}^n\mathfrak u_i(t)^2e_i^-
  &=\int_0^t\sum_{i\in\A(s)}
       \varrho_h\mathfrak u_i(t)^2\theta_i^-(s)^{\varsigma_h}
       q_i(s)\,\d s\\
 &\le\int_0^t\sum_{i\in\A(s)}
       \varrho_h\mathfrak u_i(t)^2\theta_i^-(s)^{\varsigma_h}
       |\bfv^i(s)-\bfs U(s)|\,\d s\\
 &\le\varrho_h R_U(t)^{1/3}
       \left(\sum_{i=1}^n\mathfrak u_i(t)^3
                 \int_0^t\theta_i^-(s)^{3\varsigma_h/2}\,\d s\right)^{2/3}.
 \end{aligned}
\]
Summing the individual volume bounds, multiplying by $\pi$, and taking the positive part of the
lower bound proves \eqref{eq:forecast_volume_envelope}.
\end{proof}

The lower volume bound uses \eqref{eq:mass_loss_excess} because the excess $p_i$ is bounded by the actual ocean--ice slip in \eqref{eq:exchange_error_bounds}. 
For the shortfall $q_i$, the same argument based on mass loss would require $q_i\le\Ueff^i$. This can fail: if $\bfv^i=\bfu^i$ and $c_i^->2U_{\min}$, then
$q_i=c_i^--U_{\min}>U_{\min}=\Ueff^i$.
The upper volume bound therefore uses $q_i\le|\bfv^i-\bfs U|$ and Lemma~\ref{lem:weighted_bounds}.


\subsection{Melting prediction along approximate trajectories}

Here we present an alternative estimate for the ice volume.
We first prescribe a common drift velocity and use it to construct an
approximate trajectory for each floe. We then evaluate both melt laws along
these paths using a reference temperature field and compute the resulting
volume from the initial floe sizes. Finally, we bound the difference
between this approximation and the actual volume. The bound accounts
for the error in the reference temperature field, the effect of position
error on the sampled temperature and ocean velocity, and the direct
effect of velocity error on the ocean--ice slip.

Let $\bfs U$ be the spatially uniform reference velocity chosen above,
and let $\widehat\Theta(t,\bfx)\ge0$ be a bounded reference temperature
field with the same $C^1$ regularity as $\Theta$. For each floe $i$,
define the approximate trajectory
\begin{equation}\label{eq:reference_paths}
 \bfy^i(t)=\bfx^i(0)+\int_0^t\bfs U(s)\,\d s.
\end{equation}
We use $\bfs U(t)$ as the floe velocity in the basal exchange speed.
If $\bfs U=0$, each approximate trajectory stays at the initial floe
center, so the reference temperature is sampled at that fixed position.

Specifically, set
\[
\widehat\Theta_i(t)=\widehat\Theta(t,\bfy^i(t)),\quad
\widehat U_{\mathrm{eff}}^i(t)=U_{\min}+|\bfu(t,\bfy^i(t))-\bfs U(t)|,\quad
\widehat\M_r^i=\varrho_r\widehat\Theta_i^{\varsigma_r},\quad
\widehat\M_h^i=\varrho_h\widehat U_{\mathrm{eff}}^i\widehat\Theta_i^{\varsigma_h}.
\]
This leads to a prediction of the remaining ice volume:
\[
\widehat V_U(t)
=\pi\sum_{i=1}^n
 \left(R_i-\int_0^t\widehat\M_r^i(s)\,\d s\right)_+^2
 \left(H_i-\int_0^t\widehat\M_h^i(s)\,\d s\right)_+.
\]

For $k\in\{r,h\}$, define
\[
 \begin{aligned}
 e_k(t)&=\varrho_k\bigl\|\Theta(t,\cdot)^{\varsigma_k}-\widehat\Theta(t,\cdot)^{\varsigma_k}\bigr\|_\infty,\\
 L_k(t)&=\varrho_k\bigl\|\nabla\bigl(\widehat\Theta(t,\cdot)^{\varsigma_k}\bigr)\bigr\|_\infty,
 \qquad L_o(t)=\|\nabla\bfu(t,\cdot)\|_\infty,\\
 P_h(t)&=\varrho_h\|\Theta(t,\cdot)\|_\infty^{\varsigma_h},\\
 c_U(t)&=U_{\min}+\sup_{\bfx\in\Omega}|\bfu(t,\bfx)-\bfs U(t)|.
 \end{aligned}
\]
All spatial norms are over $\Omega$.
When the actual temperature is uncertain, prescribed upper bounds for $e_k$ and $P_h$ may be used in all subsequent estimates. In particular, the bounds $\underline\Theta\le\Theta\le\overline\Theta$ give
\[
 \begin{aligned}
 e_k(t)&\le\varrho_k\left\|
   \max\left\{
    \left|\underline\Theta(t,\cdot)^{\varsigma_k}-\widehat\Theta(t,\cdot)^{\varsigma_k}\right|,
    \left|\overline\Theta(t,\cdot)^{\varsigma_k}-\widehat\Theta(t,\cdot)^{\varsigma_k}\right|
   \right\}\right\|_\infty,\qquad k\in\{r,h\},\\
 P_h(t)&\le\varrho_h\|\overline\Theta(t,\cdot)\|_\infty^{\varsigma_h}.
 \end{aligned}
\]
Set
\[
 E_\Theta(T)=\sum_{i=1}^n\int_0^T[2R_iH_ie_r(t)+R_i^2c_U(t)e_h(t)]\,\d t
\]
and
\begin{equation}\label{eq:forecast_weights}
 K_i(s)=R_i^2P_h(s)+\int_s^T
       [2R_iH_iL_r(t)+R_i^2(c_U(t)L_h(t)+P_h(t)L_o(t))]\,\d t.
\end{equation}
The first term in $K_i$ measures the immediate effect of velocity error on basal slip. The integral accounts
for its later effect on position and on the sampled temperature and ocean velocity.

Define the error bound for the volume prediction along approximate floe trajectories by
\[
 E_U^{\mathrm{fore}}(T)
 :=\pi\left[
 E_\Theta(T)+R_U(T)^{1/3}
 \left(\sum_{i=1}^n\int_0^T K_i(s)^{3/2}\,\d s\right)^{2/3}
 \right],
\]
and set
\[
 V_2^l(T):=\widehat V_U(T)-E_U^{\mathrm{fore}}(T),\qquad
 V_2^u(T):=\widehat V_U(T)+E_U^{\mathrm{fore}}(T).
\]

\begin{theorem}[Error bound for the predicted ice volume]\label{thm:reference_forecast}
Under the preceding assumptions, the predicted ice volume \(\widehat V_U(T)\) satisfies satisfies
\begin{equation}\label{eq:reference_forecast_error}
 |V(T)-\widehat V_U(T)|\le E_U^{\mathrm{fore}}(T).
\end{equation}
Consequently, $V_2^l(T)\le V(T)\le V_2^u(T)$.
These bounds remain valid when the radius or thickness vanishes in either the interacting solution or the approximate solution.
\end{theorem}
\begin{proof}
We defer the full proof to Appendix~\ref{app:forecast_error_proof}.
For each floe, we first compare the melt rates along the actual and
approximate trajectories and integrate their difference up to the earlier removal.
After that removal, one volume is zero and the other is nonincreasing,
so their difference cannot grow in absolute value. Summing over the floes
and applying Lemma~\ref{lem:weighted_bounds} with the weights
$K_i$ from \eqref{eq:forecast_weights} gives
\eqref{eq:reference_forecast_error}.
\end{proof}

We now combine the two volume intervals with the thermal upper bound
\eqref{eq:thermal_volume_upper}. Set
\[
 \underline V(T)
 :=\max\{V_1^l(T),V_2^l(T)\},\qquad
 \overline V(T)
 :=\min\{V_1^u(T),V_2^u(T),V_{\mathrm{th}}(T)\}.
\]
Since the ice volume is nonnegative, we may replace each $V_i^l$ with its positive part $(V_i^l)_+$.

\begin{corollary}\label{cor:combined_volume_interval}
Under the assumptions of Theorems~\ref{thm:volume_envelope}
and~\ref{thm:reference_forecast}, we have
\begin{equation}\label{eq:combined_volume_interval}
 \underline V(T)\le V(T)\le\overline V(T).
\end{equation}
\end{corollary}
\begin{proof}
Theorems~\ref{thm:volume_envelope} and~\ref{thm:reference_forecast}
give the two lower bounds and the two upper bounds.
Equation~\eqref{eq:thermal_volume_upper} gives the additional upper
bound $V_{\mathrm{th}}(T)$. Taking the maximum of the lower bounds
and the minimum of the upper bounds proves the claim.
\end{proof}

\begin{remark}[No uniformly sharpest bound]
\label{rem:one_floe_comparison}
We give three choices of admissible data for a dimensionless example with one floe. 
Each candidate lower bound in Corollary~\ref{cor:combined_volume_interval} is strictly the largest in at least one case, and each candidate upper bound is strictly the smallest in at least one case.
Consider the case with
\[
 T=1,\qquad R_1=2,\qquad H_1=4,\qquad
 U_{\min}=\varrho_r=\varrho_h=\varsigma_r=\varsigma_h=1.
\]
Let $\Theta=\widehat\Theta\equiv1$, $\bfs U=\bfg=0$, and take the constant fields
$\bfu=(a,0)^T$, $\bfw=(-a,0)^T$, with $a\ge0$ and $\alpha_o=\alpha_w$.
The two drag forces at the reference velocity cancel, so
\[
 R_0(1)=\frac{3m^1(0)|\bfv^1(0)|^2}{2\alpha}.
\]
Write $\eta:=R_0(1)^{1/3}$. Any specified $\eta\ge0$ can therefore be obtained by choosing
the initial speed. 
For the exact temperature bounds, define
\[
 \lambda:=
 \left(\int_0^1
   \left(1+\frac{3\rho_{ice}\pi}{2\alpha}(2-s)^2\right)^{-1/2}
   \,\d s\right)^{2/3}.
\]
Since the integrand lies strictly between zero and one, $0<\lambda<1$.
With the exact temperature bounds $\theta_1^-=\theta_1^+=1$, we have
$\mathfrak l_1(1)=\mathfrak u_1(1)=1$, $c_1^-=c_1^+=1+a$, and
\[
 \begin{aligned}
 \frac{V_1^l(1)}{\pi}&=(3-a-\lambda\eta)_+,
 &\qquad \frac{V_1^u(1)}{\pi}&=(3-a)_++\eta,\\
 \frac{V_2^l(1)}{\pi}&=(3-a)_+-4\eta,
 &\qquad \frac{V_2^u(1)}{\pi}&=(3-a)_++4\eta,\\
 \frac{V_{\mathrm{th}}(1)}{\pi}&=3.
 \end{aligned}
\]
Indeed, $\widehat V_0(1)=\pi(3-a)_+$ and
$E_0^{\mathrm{fore}}(1)=4\pi\eta$.
Here $e_r=e_h=L_r=L_h=L_o=0$, $P_h=1$, and $K_1=R_1^2P_h=4$.
For $a<3$, the robust upper bound is smaller than the thermal upper bound precisely when
$\eta<a$; the ordering is reversed when $\eta>a$.
The following choices illustrate each case in Corollary~\ref{cor:combined_volume_interval}.
\begin{enumerate}
\item If $a=1$ and $\eta=1/4$, then $[V_1^l(1),V_1^u(1)]=[\pi(2-\lambda/4),9\pi/4]$ and $[V_2^l(1),V_2^u(1)]=[\pi,3\pi]$.
The thermal upper bound is $V_{\mathrm{th}}(1)=3\pi$.
Thus $V_1^l(1)$ and $V_1^u(1)$ give both endpoints of the combined interval.
\item If $a=0$ and $\eta=1/4$, then
$[V_1^l(1),V_1^u(1)]=[\pi(3-\lambda/4),13\pi/4]$
and $[V_2^l(1),V_2^u(1)]=[2\pi,4\pi]$.
The thermal upper bound is $V_{\mathrm{th}}(1)=3\pi$.
The combined interval is $[V_1^l(1),V_{\mathrm{th}}(1)]$.
\item If $a=1$ and $\eta=0$, keep $\Theta=\widehat\Theta\equiv1$ but use the wider supplied
bounds $\underline\Theta=1/2$ and $\overline\Theta=3/2$.
Then $\mathfrak l_1(1)=1/2$, $\mathfrak u_1(1)=3/2$, and substitution gives
\[
 V_1^l(1)=\frac{\pi}{4},\qquad
 V_1^u(1)=\frac{27\pi}{4},\qquad
 V_{\mathrm{th}}(1)=\frac{63\pi}{8}.
\]
Here $\widehat V_0(1)=2\pi$ and $E_0^{\mathrm{fore}}(1)=0$, so
$V_2^l(1)=V_2^u(1)=V(1)=2\pi$. Thus the estimate along the approximate floe trajectory gives
the exact volume $V(1)=2\pi$.
The prediction $[V_2^l(1),V_2^u(1)]$ uses the exact temperature field, whereas $[V_1^l(1),V_1^u(1)]$ uses only the wider supplied bounds.
This example shows the effect of supplying different temperature information to the two estimates.
\end{enumerate}
\end{remark}

We close this section with a remark on the established ice volume bounds.
Theorem~\ref{thm:volume_envelope} establishes bounds on the ice loss from environmental bounds over admissible regions containing the floe centers.
Theorem~\ref{thm:reference_forecast} instead bounds the error in the ice volume predicted along approximate floe trajectories.
Corollary~\ref{cor:combined_volume_interval}
combines the two volume intervals with the thermal upper bound
$V_{\mathrm{th}}$.

An application of these methods requires initial floe states, calibrated coefficients, a spatial temperature forecast with error bounds, and bounds on ocean, wind, and slope forcing. 
The resulting intervals are deterministic bounds for our reduced model.
Applying the model to observed floes requires calibration of the basal melt coefficient $\varrho_h$ and the background exchange speed $U_{\mathrm{min}}$. 
Assessing the accuracy of the predicted ice volume also requires accounting for processes omitted from the model, including surface melt, area-dependent drag, tangential friction, fracture, and feedback on the ocean heat reservoir.

\section{Numerical experiments}\label{sec:num}

In this section, we consider two numerical experiments. These experiments aim to demonstrate the volume bounds developed in Section \ref{sec:prediction}.
In Example 1, the ocean and wind velocities coincide, and the temperature is constant.
We use the common velocity as the reference velocity $\bfs U$, and compare the resulting bounds.

In Example 2, we vary the environmental fields and examine how the choice of reference velocity affects the tightness of the bounds.
We select the velocity that minimizes $V_1^u$ using the grid search, using it as a constant reference velocity on $[0,T]$.

The simulations use a second-order Adams--Bashforth method for the melt equations and a velocity-Verlet-based update for the mechanical equations, with different time steps for stability and computational feasibility.

\subsection{Example 1: Constant ocean and wind velocity}

We solve the model in the domain $\Omega=[-20,20]^2$ (in meters) with double periodic boundary conditions for simulation feasibility. 
The analytical model is posed on $\mathbb R^2$; applying its estimates to this periodic computation requires lifted trajectories and a consistent periodic contact law.
For each pair of floes, the nearest periodic image remains unique throughout contact, so the contact normal used in the contact work identity is well defined.
To verify the volume bounds established in Section \ref{sec:prediction} and given the relatively small size, we suppose that $\bfg^i=\mathbf{0}$. Furthermore, to simplify the simulation, we fix $$\bfu(x,y)=\bfw(x,y)=(0.3,0)^T,$$
with $[\bfu]=[\bfw]=\text{ms}^{-1}$. We generate $n=100$ floes, for which the radii satisfy
$$R_i\in [R_{\min},R_{\max} ],\quad R_{\min}=0.5\text{m},\quad R_{\max} = 2\text{m},$$
and are sampled from a power-law distribution with exponent 2. 
All floes initially have thickness $2\,\mathrm m$.
The ice density is fixed as $917\text{kgm}^{-3}$, and the floes are initialized such that their centers are located on the lattice points of an evenly spaced grid on $\Omega$. 
Additionally, non-overlap is checked using periodic distances.  A $10\times10$ grid of spacing $4\,\mathrm m$ is sufficient for radii at most $2\,\mathrm m$, provided duplicate periodic endpoints are excluded.
The model holds the ice at $T_f=-1.9^\circ\mathrm C$ and it has no independent floe-temperature state.
For this experiment, we set $U_{\min}=0.01\,\mathrm{m\,s^{-1}}$.

The constants $\varrho_r,\,\varrho_h,\,\varsigma_r,\,\varsigma_h$ are set as above. Ocean temperature is fixed at $2^\circ $C. The other constants are set as follows: $\alpha_o=5.5\times 10^{-3}\, ,\alpha_w=1.5\times 10^{-5}\,,\kappa=10^5\,,\gamma=10^4\,,f_E=0$.

\begin{figure}[H]
    \centering
    \includegraphics[width=\textwidth]{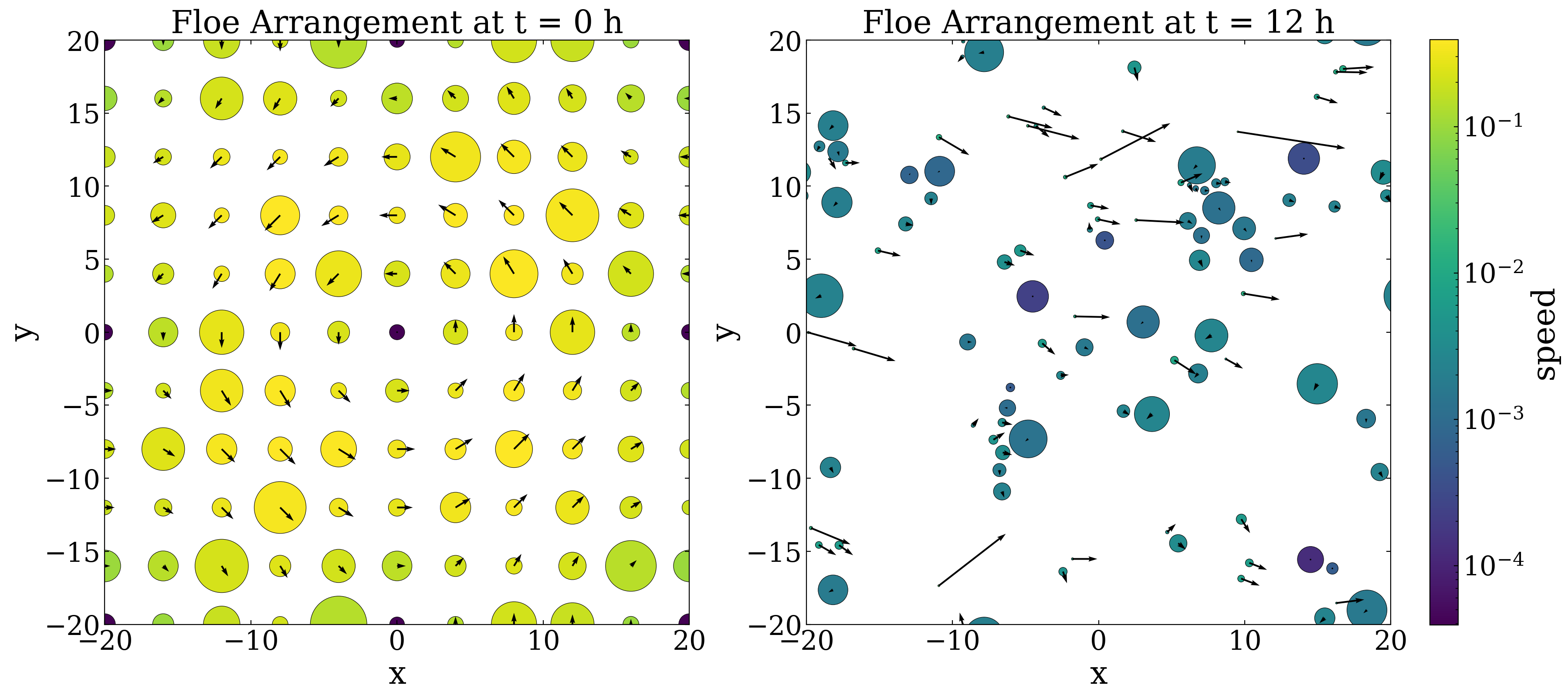}
    \caption{Initial eddy velocity field and floe positions after 12 hours. The floes drift toward the ocean velocity $\bfu$. 
    Colors indicate drift speed on the same scale in both panels.
    Arrows are scaled separately in each panel so that drift directions
    remain visible at low speeds. }
    \label{fig:sim_1_ice_floe_movement}
\end{figure}

We use a prescribed eddy as the initial velocity of the floes. Following example 1 in \cite{deng2026particle}, we define the initial floe velocities as
$$\bfv^i(0)=0.4\mathbf{e}(\hat{\bfx}^i),\quad \mathbf{e}(\hat{\bfx}^i)=\varphi(\hat{\bfx}^i)\begin{pmatrix}
    -\sin(\hat{y}^i)\\
    \sin(\hat{x}^i)
\end{pmatrix},$$
where
$$\varphi(\hat{\bfx}^i)=\exp\left(\frac{\cos (\hat{x}^i)+\cos(\hat{y}^i)-2}{4.41}\right),\quad \hat{x}^i=\frac{\pi x^i}{x_{\max}},\quad \hat{y}^i=\frac{\pi y^i}{y_{\max}}.$$
We run the simulation over two days (Figure~\ref{fig:sim_1_ice_floe_movement}).

To compute the bounds, we choose $\bfu$ as our reference velocity $\bfs U$. Then, we have
$$\bfs B_U(t,\bfx,m)=\mathbf{0},\qquad b_i(t)=0,\qquad c_i^-(s)=c_i^+(s)=U_{\min}.$$
The temperature $T_o$ is constant, and we simply take 
\[
\underline{\Theta}=\widehat\Theta=\overline\Theta=\Theta.
\]
Along the approximate trajectories, the ocean--ice slip is zero.
The predicted volume therefore coincides with the thermal upper bound:
\[
 \widehat V_U(t)=V_{\mathrm{th}}(t)
 =\pi\sum_{i=1}^n
   \bigl(R_i-\varrho_r\Theta^{\varsigma_r}t\bigr)_+^2
   \bigl(H_i-\varrho_h U_{\min}\Theta^{\varsigma_h}t\bigr)_+.
\]
Moreover, \eqref{eq:forecast_weights} reduces to
$K_i=R_i^2\varrho_h\Theta^{\varsigma_h}$, and the volume error bound becomes
\begin{equation}\label{eq:sim_1_fore}
 E_U^{\mathrm{fore}}(T)
 =\pi\varrho_h\Theta^{\varsigma_h}
  \left(\frac{3E_U(0)}{\alpha}\right)^{1/3}
  \left(\sum_{i=1}^n R_i^3\right)^{2/3}T^{2/3}.
\end{equation}

\begin{figure}[H]
    \centering
    \includegraphics[width=0.8\textwidth]{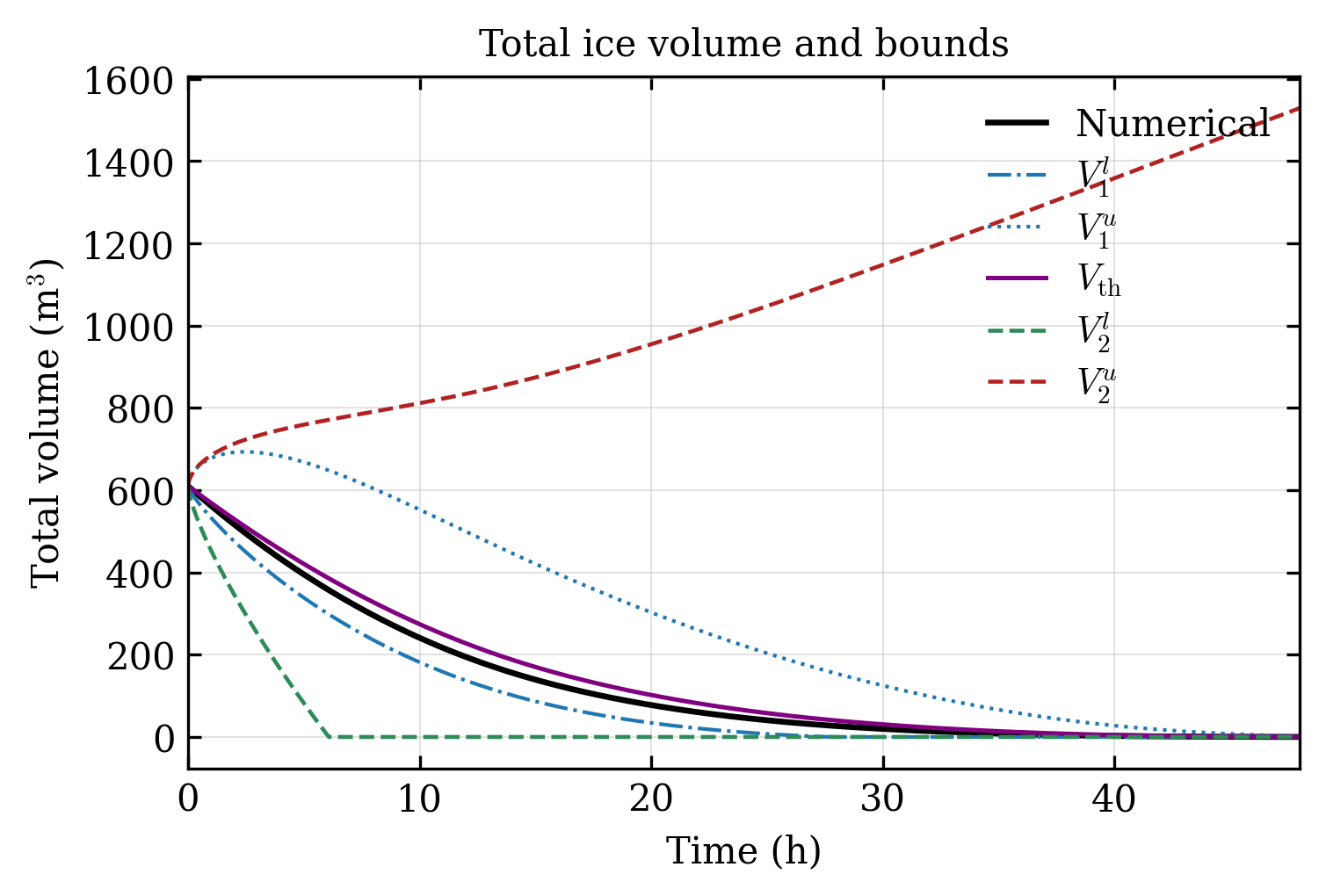}
    \caption{Simulated total ice volume (black curve) in Example~1 and the bounds
from Section~\ref{sec:prediction}. 
Among the bounds shown, $V_1^l$ and $V_{\mathrm{th}}$ give the tightest lower and upper bounds, respectively.}
    \label{fig:sim_1}
\end{figure}

To compare the volume bounds, we report their $L^2$ deviations from the computed volume, normalized by the initial volume: $\|(\tilde{V}-V)/V(0)\|_{L^2(0,T)}$, where $\tilde{V}$ is a volume bound. 
The results are summarized in the Table~\ref{table:sim_1}.

\begin{table}[htbp]
\centering
\begin{tabular}{cccc}
\toprule
\multicolumn{2}{c}{Bound} & $\|(\tilde{V}-V)/V(0)\|_{L^2}$ & $\max |(\tilde{V}-V)/V(0)|$ \\
\midrule
\multirow{3}{*}{Upper}&$V_\mathrm{th}$ & 13.50 & 5.46\% \\
&$V^u_1$ & 129.85 & 51.12\% \\
&$V_2^u$ & 689.61 & 250.55\% \\
\midrule
\multirow{2}{*}{Lower}&$V^l_1$ & 25.42 & 9.85\% \\
& $V^l_2$ & 95.22 & 58.49\% \\
\bottomrule
\end{tabular}
\caption{$L^2$ and maximum absolute error for volume bounds, normalized by initial volume $V_0$.}
\label{table:sim_1}
\end{table}

As shown by (\ref{eq:sim_1_fore}) and Figure \ref{fig:sim_1}, $E_U^{\mathrm{fore}}$ increases in time, leading to a significant error compared to the other bounds. 
Overall, the sharpest bounds are given by $$\underline{V}(T)=V_1^l(T),\quad \overline{V}(T)=V_\mathrm{th}(T).$$

\begin{figure}[H]
    \centering
    \includegraphics[width=0.85\linewidth]{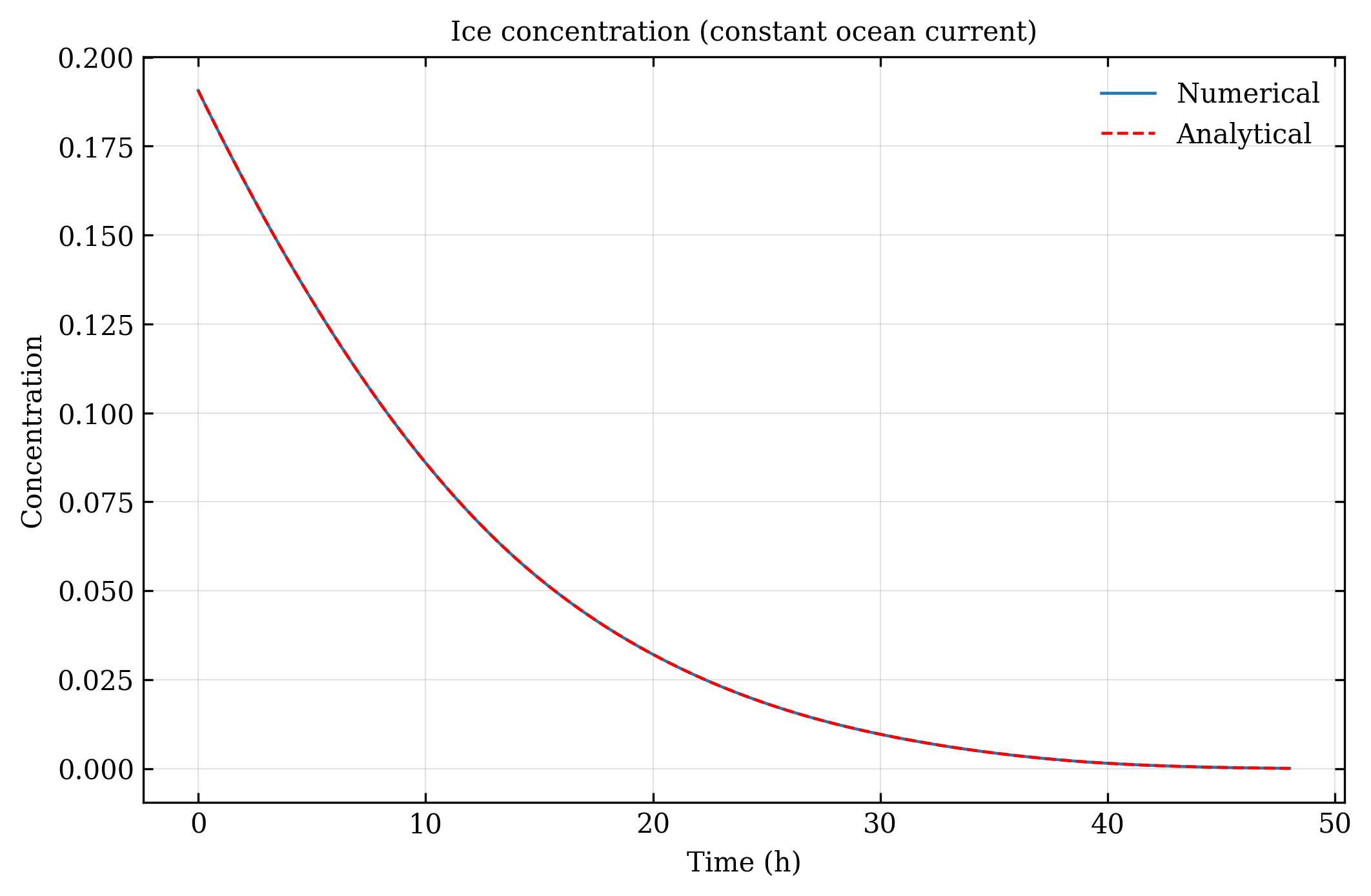}
    \caption{Smooth decay of ice concentration.}
    \label{fig:sim_1_conc}
\end{figure}

Furthermore, since the radius melting is independent of forcing dynamics, for this specific experiment, for each active floe, the radius is 
\[
 r^i(t)=R_i-\varrho_r\Theta^{\varsigma_r}t,\qquad t<\tau^i.
\]
The sea ice concentration, defined as normalized sum of projected floe areas, is therefore
\[
 \frac{A(t)}{|\Omega|}=\frac{\pi}{|\Omega|}
 \sum_{i:\,t<\tau^i}\left(R_i-\varrho_r\Theta^{\varsigma_r}t\right)_+^2.
\]
See Figure~\ref{fig:sim_1_conc} for the evolution of the concentration. 
The active-set restriction is essential if thickness vanishes first. This quantity decreases continuously between removals and has no jump at radius extinction; thickness extinction can produce a downward jump. As floes may overlap, it need not equal the fractional area of the union of the floe disks. We assume the overlap is negligible here for the calculation of ice concentration.

\subsection{Example 2: Constant ocean and wind in different directions}

In Example 1, $V_{\mathrm{th}}$ provides the tightest upper bound. 
Keeping the other settings of Example~1, we change the wind velocity, temperature field, and the initial velocity of floes.
Specifically, we set
\[
 \bfw(x,y)=(-10,0)^T,\quad
 T_o(t)=2+0.5\cos\left(\frac{2\pi t}{86400}\right),\quad
 \bfv^i(0)=\boldsymbol{u}_*+0.01\mathbf{e}(\hat{\bfx}^i),
\]
where
\[
 \boldsymbol{u}_*=
 \frac{\sqrt{\alpha_o}\bfu+\sqrt{\alpha_w}\bfw}
 {\sqrt{\alpha_o}+\sqrt{\alpha_w}}.
\]

We consider three reference velocities when computing the
volume bounds: $U\in\{\bfu,\bfw,\bfs U_\ast\}$.
The reference velocity $\bfs U_\ast$ is chosen as a constant velocity
\[
 \bfs U_\ast(t)\in\operatorname*{arg\,min}_{\mathbf a\in\mathcal U}
 V_1^u(t;\bfs U\equiv \mathbf a),
\]
where the search is performed over the following grid:
\[
 \mathcal U=\mathcal X\times\mathcal Y,\quad
 \mathcal X=\left\{-12+\frac{7n}{8}:0\le n\le16\right\},\quad
 \mathcal Y=\left\{-3+\frac{3n}{5}:0\le n\le10\right\}.
\]
Herein, for each evaluation time $t\in [0,T]$, every candidate is used as a constant reference velocity over the entire interval $[0,t]$. 
The pointwise minimum of these valid upper bounds remains an upper bound. 
A sequence of timewise minimizers must not be substituted as a single reference velocity without the required $C^1$ regularity and acceleration residual. Ties may be broken deterministically.

On the other hand, by our choice of $T_o$, it follows that $\Theta(t,\bfx)$ is in fact independent of $\bfx$. We simply set the bounds to be $\Theta$, i.e.,
$$\underline\Theta(t)=\Theta(t)=\overline\Theta(t).$$
The fact that temperature is constant across space also means that we have
$$\theta_i^-(t)=\theta_i^+(t)=\Theta(t).$$
With $\widehat\Theta=\Theta$, the temperature-field error vanishes. However, for a constant reference velocity $\bfs U$,
\[
 \widehat V_U(t)=\pi\sum_{i=1}^n
 \left(R_i-\varrho_r\int_0^t\Theta(s)^{\varsigma_r}\,\d s\right)_+^2
 \left(H_i-\varrho_h(U_{\min}+|\bfu-\bfs U|)
                   \int_0^t\Theta(s)^{\varsigma_h}\,\d s\right)_+.
\]
Consequently $\widehat V_U(t)\le V_{\mathrm{th}}(t)$, with equality when $\bfs U=\bfu$.
\begin{table}[H]
\centering
\begin{tabular}{cc c cc}
\toprule
\multicolumn{2}{c}{Bound} & $U$ & {$\|(\tilde{V}-V)/V(0)\|_{L^2}$} & {$\max |(\tilde{V}-V)/V(0)|$} \\
\midrule
\multirow{7}{*}{Upper}
& $V_\mathrm{th}$ & -- & 26.92 & 10.56\% \\
\cmidrule(lr){2-5}
& \multirow{3}{*}{$V^u_1$}
& $\bfu$ & 168.71 & 64.16\% \\
& & $\bfw$ & 1177.12 & 434.79\% \\
& & $\bfs U_\ast$ & 3.70 & 1.52\% \\
\cmidrule(lr){2-5}
& \multirow{3}{*}{$V^u_2$}
& $\bfu$ & 820.57 & 301.03\% \\
& & $\bfw$ & 8370.58 & 3391.77\% \\
& & $\bfs U_\ast$ & 4.87 & 1.68\% \\
\midrule
\multirow{6}{*}{Lower}
& \multirow{3}{*}{$V^l_1$}
& $\bfu$ & 20.24 & 7.80\% \\
& & $\bfw$ & 121.04 & 90.21\% \\
& & $\bfs U_\ast$ & 3.20 & 1.33\% \\
\cmidrule(lr){2-5}
& \multirow{3}{*}{$V^l_2$}
& $\bfu$ & 53.80 & 31.04\% \\
& & $\bfw$ & 121.33 & 90.69\% \\
& & $\bfs U_\ast$ & 2.916 & 1.26\% \\
\bottomrule
\end{tabular}
\caption{$L^2$ and maximum absolute error for volume bounds, normalized by initial volume $V_0$.}
\label{table:sim_2}
\end{table}

The reported results suggest that choosing the reference velocity carefully can substantially tighten the bounds.

\section{Conclusion}\label{sec:conclusion}

This paper studies a floe model that couples translation, lateral recession, and basal thinning under marine heatwaves.
Our model has a unique local solution under regularity conditions on the environmental fields and a non-grazing condition on initial contacts.
This local solution can be continued uniquely to a global solution through removals under non-rafting and non-grazing conditions, with the latter also imposed at limiting contacts of surviving floes.
When $\gamma=0$, the contact force is locally Lipschitz across contact onset, so no non-grazing contact assumption is needed in either case.
Here, the pointwise velocity bound keeps the forces on surviving floes bounded, so their full states have finite limits at extinction.

Next, we establish estimates for the relative energy and integrated relative speed in terms of initial states and environmental fields.
For any given solution obeying the removal rule, we provide two types of ice volume bounds using these estimates.
The first estimate predicts the ice volume by combining bounds on the environmental fields over admissible regions containing the floe centers with an integral bound on relative speed.
The second estimate predicts the ice volume by evaluating the melt laws along approximate floe trajectories and then bounding the error in the predicted volume.

Numerical experiments validate these established volume bounds and show how tight they are.
In Example~2, we select a constant reference velocity at each evaluation
horizon by minimizing the regional upper bound $V_1^u$ over a finite
grid. Compared with using the ocean velocity $\bfu$ as the reference,
this choice reduces the maximum gap between $V_1^u$ and the simulated
volume from $64.16\%$ to $1.52\%$ of the initial volume $V(0)$.
The same choice also tightens the regional lower bound $V_1^l$:
its maximum gap from the simulated volume decreases from $7.80\%$
to $1.33\%$ of $V(0)$ (Table~\ref{table:sim_2}).

For future work, one direction is to extend the model to more realistic floe configurations by incorporating rotation, irregular shapes, and coupled melting and refreezing. 
These features alter contact geometry, force and torque transmission, and the evolving floe size distribution, with consequences for collective motion and deformation. 
Another goal is to derive effective sea ice rheologies from these interactions and determine how thermodynamic changes influence the transition between fragmented floe dynamics and a mechanically connected ice cover, leading to more accurate ice predictions in the polar region.

\appendix
\counterwithin{table}{section}

\section{Temperature-index melt closure}\label{app:melt}

We relate the melt rates to lateral recession and a balance of heat flux from the ocean to the ice, and record the reference
parameter values. We prescribe $T_o\ge T_f$, keep each floe internally at $T_f$, and omit
freezing, internal heat storage, surface ablation and feedback on the
prescribed ocean heat reservoir.  The two interface recession speeds in
\eqref{eq:meltrh} therefore describe a reduced closure over the duration of the event,
not a complete ice--ocean energy balance.

\paragraph{Lateral recession.}
For a circular floe undergoing uniform edge recession, the
Maykut--Perovich parameterization given in 
\citep[Eq.~(2)]{bateson2020impact}
yields
\[
 -\dot r^i=\M_r^i=\varrho_r\Theta_i^{\varsigma_r},
 \qquad
 \varrho_r=1.6\times10^{-6}\,
   \mathrm{m\,s^{-1}}({}^{\circ}\mathrm{C})^{-1.36},
 \qquad
 \varsigma_r=1.36.
\]

\paragraph{Basal recession.}
Let $F_b^i$ denote the heat flux from the ocean to the base of floe~$i$,
taken to be positive when heat enters the ice.
Equation~(2) of \cite{maykut1995solar} gives the bulk heat flux
\[
 F_b^i=\rho_o c_w c_H u_*^i\Theta_i,
\]
where $\rho_o$ and $c_w$ are the density and specific heat capacity of
seawater, $c_H$ is the heat transfer coefficient, and $u_*^i$ is the
friction velocity at the ice--ocean interface.
Since we assume that the ice is uniformly at $T_f$, the conductive flux
in the basal heat balance of \cite[Eq.~(3)]{maykut1995solar} vanishes.
The balance therefore reduces to
\[
 \rho_{ice}L_i\M_h^i=F_b^i,
\]
where $\rho_{ice}$ is the ice density and $L_i$ is the specific latent
heat of fusion, for which we use the fresh-ice value.

\cite[Eq.~(1)]{maykut1995solar} estimate the friction velocity from
ice drift relative to the surface geostrophic current using a complex
geostrophic drag relation. For the resolved ice--ocean relative
velocity in our model, we use the approximation
$\sqrt{C_D^{io}}\,|\bfu^i-\bfv^i|$ of \cite[Eq.~(4)]{nihashi2008bulk}, where $C_D^{io}$ is the ice--ocean drag coefficient.
In the CICE simulations of Bateson et al., the heat available from the upper ocean limits the combined lateral and basal melt \citep[Sec.~2.1]{bateson2020impact}.
Gu et al.\ use a minimum ocean friction velocity of $5\times10^{-4}\,\mathrm{m\,s^{-1}}$ in their Icepack simulations
\citep[Table~1]{gu2022landfast}.
We adopt this value for $u_{*,\min}$ and add it to the resolved friction velocity to represent unresolved turbulent exchange.
Thus,
\[
 u_*^i
 =u_{*,\min}+\sqrt{C_D^{io}}\,|\bfu^i-\bfv^i|
 =\sqrt{C_D^{io}}\,\Ueff^i,
 \qquad
 U_{\min}=\frac{u_{*,\min}}{\sqrt{C_D^{io}}}.
\]
Substituting into the basal heat balance gives
\[
 -\dot h^i=\M_h^i=\varrho_h\Ueff^i\Theta_i,
 \qquad \varsigma_h=1.
\]

\begin{table}[H]
\centering
\caption{Reference values used in the basal melt law.}
\label{tab:basal_melt_parameters}
\begin{tabular}{@{}lll@{}}
\hline
Symbol & Value & Source \\
\hline
$\rho_o$
  & $1025\,\mathrm{kg\,m^{-3}}$
  & \citep[Sec.~3f]{parkinson1979large} \\
$c_w$
  & $3980\,\mathrm{J\,kg^{-1}}({}^{\circ}\mathrm{C})^{-1}$
  & \citep[p.~24698]{maykut1995solar} \\
$c_H$
  & $0.0057$
  & \citep[para.~7]{mcphee2003ocean} \\
$C_D^{io}$
  & $5.5\times10^{-3}$
  & \citep[Table~1]{hunke1997elastic} \\
$\rho_{ice}$
  & $917\,\mathrm{kg\,m^{-3}}$
  & \citep[Table~1]{huwald2005multilayer} \\
$L_i$
  & $3.34\times10^5\,\mathrm{J\,kg^{-1}}$
  & \citep[Table~1]{huwald2005multilayer} \\
$u_{*,\min}$
  & $5\times10^{-4}\,\mathrm{m\,s^{-1}}$
  & Minimum value in \citep[Table~1]{gu2022landfast} \\
\hline
\end{tabular}
\end{table}

The values in Table~\ref{tab:basal_melt_parameters} give
\[
\begin{aligned}
 \varrho_h
 &=\frac{\rho_o c_w c_H\sqrt{C_D^{io}}}{\rho_{ice}L_i}
  =\frac{1025(3980)(0.0057)\sqrt{5.5\times10^{-3}}}
         {917(3.34\times10^5)}
  \simeq 5.63\times10^{-6}\,({}^{\circ}\mathrm{C})^{-1},\\
 U_{\min}
 &=\frac{u_{*,\min}}{\sqrt{C_D^{io}}}
  =\frac{5\times10^{-4}}{\sqrt{5.5\times10^{-3}}}
  \simeq 6.74\times10^{-3}\,\mathrm{m\,s^{-1}}.
\end{aligned}
\]

\section{Mechanical estimate}\label{app:motion_budget_explicit}

We prove Proposition~\ref{prop:motion_budget} by computing the relative work and integrating the resulting differential inequality.

\begin{proof}[Proof of Proposition~\ref{prop:motion_budget}]
Write $\bfs e^i=\bfv^i-\bfs U$ and $\bfz_o^i=\bfu^i-\bfs U$, $\bfz_w^i=\bfw^i-\bfs U$.

\medskip
\noindent\textbf{$\bullet$ Step 1} (Identity for the relative energy).
On an interval with a fixed active set, differentiating the kinetic part of \eqref{eq:relative_energy}
and using \eqref{eq:dem} gives
\begin{equation}\label{eq:relative_kinetic_derivative}
 \begin{aligned}
 \frac{\d}{\d t}\left(\frac12\sum_i m^i|\bfs e^i|^2\right)
 ={}&\sum_i\left(\sum_{j\ne i}\bff^{ij}+\bfF_e^i-m^i\dot{\bfs U}\right)
               \cdot\bfs e^i\\
 &+\frac12\sum_i\dot m^i|\bfs e^i|^2.
 \end{aligned}
\end{equation}
Antisymmetry and the gap equation give the identity for the work of the contact forces and the derivative of the spring energy, respectively:
\[
 \begin{aligned}
 \sum_i\sum_{j\ne i}\bff^{ij}\cdot\bfs e^i
 &=\sum_{i<j}\chi^{ij}(\kappa\delta^{ij}-\gamma q^{ij})q^{ij},\\
 \frac{\d}{\d t}\left(\frac\kappa2\sum_{i<j}(-\delta^{ij})_+^2\right)
 &=\kappa\sum_{i<j}\chi^{ij}\delta^{ij}
                  (-q^{ij}+\M_r^i+\M_r^j).
 \end{aligned}
\]
Substituting the identity for the work of the contact forces into \eqref{eq:relative_kinetic_derivative} and adding the derivative of the spring energy yields
\begin{equation}\label{eq:relative_energy_identity}
 \begin{aligned}
 \dot E_U={}&-\gamma\sum_{i<j}\chi^{ij}(q^{ij})^2
       +\frac12\sum_i\dot m^i|\bfs e^i|^2
       +\kappa\sum_{i<j}\chi^{ij}\delta^{ij}(\M_r^i+\M_r^j)\\
 &+\sum_i(\bfF_e^i-m^i\dot{\bfs U})\cdot\bfs e^i.
 \end{aligned}
\end{equation}
The environmental work is exactly
\begin{equation}\label{eq:environmental_work_identity}
 \begin{aligned}
 (\bfF_e^i-m^i\dot{\bfs U})\cdot\bfs e^i
 ={}&\sum_{\nu\in\{o,w\}}\alpha_\nu \Phi(\bfz_\nu^i-\bfs e^i)\cdot\bfs e^i
     +m^i(\bfg^i-\dot{\bfs U}-f_EJ\bfs U)\cdot\bfs e^i\\
 ={}&\sum_{\nu\in\{o,w\}}\alpha_\nu
       [\Phi(\bfz_\nu^i-\bfs e^i)-\Phi(\bfz_\nu^i)]\cdot\bfs e^i
     +\bfs B_U(t,\bfx^i,m^i)\cdot\bfs e^i.
 \end{aligned}
\end{equation}
Here $J\bfv^i=J\bfs e^i+J\bfs U$ and $J\bfs e^i\cdot\bfs e^i=0$.
For $\Phi(\bfz)=|\bfz|\bfz$,
\[
 (\Phi(\bfx)-\Phi(\bfy))\cdot(\bfx-\bfy)
 =\frac{|\bfx|+|\bfy|}{2}
       \bigl(|\bfx-\bfy|^2+(|\bfx|-|\bfy|)^2\bigr)
 \ge\tfrac12|\bfx-\bfy|^3.
\]
Taking the two arguments to be $\bfz_\nu^i$ and $\bfz_\nu^i-\bfs e^i$, and using the definition of $b_i$, gives
\begin{equation}\label{eq:environmental_work_bound}
 (\bfF_e^i-m^i\dot{\bfs U})\cdot\bfs e^i
 \le-\frac\alpha2|\bfs e^i|^3
       +\bfs B_U(t,\bfx^i,m^i)\cdot\bfs e^i
 \le-\frac\alpha2|\bfs e^i|^3+b_i|\bfs e^i|.
\end{equation}
The contact terms in \eqref{eq:relative_energy_identity} are nonpositive.
Moving the mass derivative term to the left and applying
\eqref{eq:environmental_work_bound}, we obtain
\[
 \dot E_U+\frac{\alpha}{2}\sum_i|\bfs e^i|^3
 +\frac12\sum_i(-\dot m^i)|\bfs e^i|^2
 \le\sum_i b_i|\bfs e^i|.
\]

\medskip
\noindent\textbf{$\bullet$ Step 2} (Forcing and integration).
For $b,z\ge0$, direct maximization of $bz-\alpha z^3/6$ gives
\[
 bz\le\frac\alpha6z^3+\frac{2\sqrt2}{3\sqrt\alpha}b^{3/2}.
\]
Consequently,
\begin{equation}\label{eq:relative_differential_inequality}
 \frac3\alpha\dot E_U
 +\sum_i\left(
     |\bfs e^i|^3+\frac{3}{2\alpha}(-\dot m^i)|\bfs e^i|^2
   \right)
 \le\frac{2\sqrt2}{\alpha^{3/2}}\sum_i b_i^{3/2}.
\end{equation}
Integration on a compact phase interval $[s,t]$ yields
\begin{equation}\label{eq:phase_relative_budget}
 \begin{aligned}
 \frac3\alpha E_U(t)
 &+\int_s^t\sum_{i\in\A(\sigma)}
   \left(
     |\bfs e^i(\sigma)|^3
     +\frac{3}{2\alpha}(-\dot m^i(\sigma))|\bfs e^i(\sigma)|^2
   \right)\,\d\sigma\\
 &\le\frac3\alpha E_U(s)
   +\frac{2\sqrt2}{\alpha^{3/2}}
      \int_s^t\sum_{i=1}^n b_i(\sigma)^{3/2}\,\d\sigma.
 \end{aligned}
\end{equation}
This estimate holds on every compact subinterval of a phase between removals.
If the phase begins at $s$ and ends at a finite removal time $\tau$, letting
$t\uparrow\tau$ gives
\[
 \int_s^\tau \sum_i |\bfs e^i(t)|^3\,\d t<\infty.
\]
Since $\tau-s<\infty$, H\"older's inequality gives
$\bfs e^i\in L^1(s,\tau)$ for each active label $i$.
The boundedness of $\bfs U$ then gives
$\bfv^i=\bfs e^i+\bfs U\in L^1(s,\tau)$.
Hence each center has a finite left limit,
\[
 \lim_{t\uparrow\tau}\bfx^i(t)
 =\bfx^i(s)+\int_s^\tau \bfv^i(t)\,\d t.
\] The differential inequality also shows that
$E_U$ minus its integrated forcing is nonincreasing and bounded below, so $E_U(\tau^-)$ exists and is finite.
For the retained states of the given solution,
\[
 \begin{aligned}
 E_U(\tau)
 &=\lim_{t\uparrow\tau}\left[
       \frac12\sum_{i\in\A(\tau)}m^i|\bfv^i-\bfs U|^2
       +\frac\kappa2\sum_{\substack{i<j\\i,j\in\A(\tau)}}(-\delta^{ij})_+^2\right]\le E_U(\tau^-).
 \end{aligned}
\]
Thus deleting nonnegative kinetic and spring contributions cannot increase energy. Taking endpoint limits and
summing \eqref{eq:phase_relative_budget} over the finitely many intervals between removals proves
\eqref{eq:weighted_speed_result}. No pointwise speed estimate is used.

\end{proof}

\section{Error of the melting laws along approximate floe trajectories}\label{app:forecast_error_proof}

We compare the melt rates during the common lifetime, use monotonicity after the earlier removal, and
apply Lemma~\ref{lem:weighted_bounds} to the resulting velocity weight.

\begin{proof}[Proof of Theorem~\ref{thm:reference_forecast}]
We first compare the melt rates, then the volumes through the earlier removal, and finally apply
the weighted estimate for accumulated relative motion.

\medskip
\noindent\textbf{$\bullet$ Step 1} (Comparison of melt rates).
Define the reference sizes by
\[
 \widehat r^i(t)=\left(R_i-\int_0^t\widehat\M_r^i(s)\,\d s\right)_+,
 \qquad
 \widehat h^i(t)=\left(H_i-\int_0^t\widehat\M_h^i(s)\,\d s\right)_+,
\]
and their removal time by
\[
 \widehat\tau^i:=\inf\{t>0:\widehat r^i(t)=0\ \text{or}\ \widehat h^i(t)=0\}.
\]
Set $\sigma_i=T\wedge\tau^i\wedge\widehat\tau^i$.
During the common active lifetime, the temperature difference is bounded by
\[
 \varrho_k|\Theta(t,\bfx^i)^{\varsigma_k}-\widehat\Theta(t,\bfy^i)^{\varsigma_k}|
 \le e_k(t)+L_k(t)|\bfx^i-\bfy^i|,\qquad k\in\{r,h\}.
\]
The difference between the exchange speeds satisfies
\[
 |\Ueff^i-\widehat U_{\mathrm{eff}}^i|
 \le L_o|\bfx^i-\bfy^i|+|\bfv^i-\bfs U|.
\]
For the basal melt rate, use the identity
\[
 \M_h^i-\widehat\M_h^i
 =\varrho_h\widehat U_{\mathrm{eff}}^i
      (\Theta_i^{\varsigma_h}-\widehat\Theta_i^{\varsigma_h})
  +\varrho_h\Theta_i^{\varsigma_h}
      (\Ueff^i-\widehat U_{\mathrm{eff}}^i).
\]
The bounds $\widehat U_{\mathrm{eff}}^i\le c_U$ and
$\varrho_h\Theta_i^{\varsigma_h}\le P_h$, together with the preceding temperature estimate, give
\begin{equation}\label{eq:forecast_melt_comparison}
 \begin{aligned}
 |\M_r^i-\widehat\M_r^i|&\le e_r+L_r|\bfx^i-\bfy^i|,\\
 |\M_h^i-\widehat\M_h^i|&\le c_Ue_h+(c_UL_h+P_hL_o)|\bfx^i-\bfy^i|
                                  +P_h|\bfv^i-\bfs U|.
 \end{aligned}
\end{equation}

\medskip
\noindent\textbf{$\bullet$ Step 2} (Volume comparison through removal).
The elementary product estimate is
\[
 |r^2h-\widetilde r^2\widetilde h|
 \le2R_iH_i|r-\widetilde r|+R_i^2|h-\widetilde h|
 \quad(0\le r,\widetilde r\le R_i,\ 0\le h,\widetilde h\le H_i).
\]
Indeed, write $r^2h-\widetilde r^2\widetilde h=(r^2-\widetilde r^2)h+
\widetilde r^2(h-\widetilde h)$ and use $|r^2-\widetilde r^2|\le2R_i|r-\widetilde r|$.
The common initial sizes and the integrated melt laws give, for $0\le t\le\sigma_i$,
\[
 \begin{aligned}
 |r^i(t)-\widehat r^i(t)|&\le\int_0^t|\M_r^i(s)-\widehat\M_r^i(s)|\,\d s,\\
 |h^i(t)-\widehat h^i(t)|&\le\int_0^t|\M_h^i(s)-\widehat\M_h^i(s)|\,\d s.
 \end{aligned}
\]
At a removal endpoint, sizes are interpreted as left limits; equivalently, one may use
$|a_+-b_+|\le|a-b|$ in the integrated size formulas. Combining these estimates at $t=\sigma_i$ yields
\[
 \begin{aligned}
 &\bigl|(r^i(\sigma_i))^2h^i(\sigma_i)
       -(\widehat r^i(\sigma_i))^2\widehat h^i(\sigma_i)\bigr|\\
 &\qquad\le\int_0^{\sigma_i}
       [2R_iH_i|\M_r^i(s)-\widehat\M_r^i(s)|
        +R_i^2|\M_h^i(s)-\widehat\M_h^i(s)|]\,\d s.
 \end{aligned}
\]
If $\sigma_i<T$, one product is zero and remains zero, whereas the other is nonnegative and nonincreasing.
Its later discrepancy is therefore no larger than at $\sigma_i$. Summing gives, through either extinction mode,
\begin{equation}\label{eq:forecast_common_lifetime}
 \frac1\pi|V(T)-\widehat V_U(T)|
 \le\sum_{i=1}^n\int_0^{\sigma_i}
       [2R_iH_i|\M_r^i-\widehat\M_r^i|+R_i^2|\M_h^i-\widehat\M_h^i|]\,\d t.
\end{equation}

\medskip
\noindent\textbf{$\bullet$ Step 3} (Weight for the velocity error).
The initial positions agree, so $|\bfx^i(t)-\bfy^i(t)|\le\int_0^t|\bfv^i(s)-\bfs U(s)|\,\d s$ before
actual removal. 
For any nonnegative function $A_i(t)$, Tonelli's theorem gives
\[
 \begin{aligned}
 &\int_0^{\sigma_i} A_i(t)\int_0^t|\bfv^i(s)-\bfs U(s)|\,\d s\,\d t\\
 &\hspace{0.5cm}=\int_0^{\sigma_i}\left(\int_s^{\sigma_i}A_i(t)\,\d t\right)
|\bfv^i(s)-\bfs U(s)|\,\d s 
 \le\int_0^{\sigma_i}\left(\int_s^T A_i(t)\,\d t\right)
 |\bfv^i(s)-\bfs U(s)|\,\d s.
 \end{aligned}
\]
Apply this with $A_i=2R_iH_iL_r+R_i^2(c_UL_h+P_hL_o)$ in
\eqref{eq:forecast_common_lifetime}, using \eqref{eq:forecast_melt_comparison}.
The direct velocity term supplies $R_i^2P_h$, giving precisely the weight \eqref{eq:forecast_weights}:
\[
 \frac1\pi|V(T)-\widehat V_U(T)|
 \le E_\Theta(T)+\sum_{i=1}^n\int_0^{\sigma_i}K_i(s)|\bfv^i(s)-\bfs U(s)|\,\d s.
\]
Since $\sigma_i\le T\wedge\tau^i$ and $K_i\ge0$, applying Lemma~\ref{lem:weighted_bounds} gives
\[
 \begin{aligned}
 &\sum_{i=1}^n\int_0^{\sigma_i}K_i(s)|\bfv^i(s)-\bfs U(s)|\,\d s\\
 &\hspace{0.5cm}\le\int_0^T\sum_{i\in\A(s)}K_i(s)|\bfv^i(s)-\bfs U(s)|\,\d s
 \le R_U(T)^{1/3}
       \left(\sum_{i=1}^n\int_0^T K_i(s)^{3/2}\,\d s\right)^{2/3}.
 \end{aligned}
\]
This proves \eqref{eq:reference_forecast_error}.
\end{proof}

\section*{Acknowledgements}
We acknowledge the following financial support. Q.D. is partially supported by the Tsinghua University Dushi program and the start-up funding from the Yau Mathematical Sciences Center, Tsinghua University. The work of S.-Y. Ha is
supported by National Research Foundation (NRF) grant funded by the Korea government (MSIT)
(RS-2025-00514472).
The funders had no role in the study design, data collection or
analysis, interpretation of the results, preparation of the
manuscript, or decision to submit it for publication.

\section*{Data and code availability}
The Python scripts used to run the two numerical experiments in
Section~\ref{sec:num}, compute the errors, and generate the figures
are available at
\url{https://github.com/yipyanyujulian/Sea-Ice-Floe-Dynamics-and-Volume-Loss-Estimates-under-Marine-Heatwaves}.

\section*{Declaration of generative AI and AI-assisted technologies
in the manuscript preparation process}
The authors used ChatGPT (OpenAI) to revise the English prose,
organize parts of the exposition, and help draft the schematic in
Figure~\ref{fig:two_volume_bounds}. The authors reviewed and edited
the resulting material and take full responsibility for the published work.

\bibliographystyle{elsarticle-num}
\bibliography{ref}

\end{document}